\documentclass[final,3p,times,authoryear, number,comma,square]{elsarticlePLF} 
\usepackage{amssymb,color}
\usepackage{a4wide}
\usepackage{amsmath,amsthm}
\usepackage{amssymb}

\numberwithin{equation}{section} 
\numberwithin{equation}{section}
\numberwithin{figure}{section}
\newcommand \bei {\begin{itemize}}
\newcommand \eei {\end{itemize}}
\newcommand \Pcal {\mathcal P} 
\newcommand \be   {\begin{equation}}
\newcommand \ee   {\end{equation}}
\newcommand \supp {{\text{supp }}} 
\newcommand \RR    {\mathbb{R}} 
\newcommand \del   {\partial}
 
\journal{\rm Journal of Computational Physics}

\begin{document}

\begin{frontmatter}
 
\title{Revisiting the method of characteristics 
via a convex hull algorithm}


\author{Philippe G. LeFloch$^1$ and Jean-Marc Mercier$^2$} 

\address{$^1$ Laboratoire Jacques-Louis Lions \& Centre National de la Recherche Scientifique, 
Universit\'e Pierre et Marie Curie, 
\\
4 Place Jussieu, 75252 Paris, France. 
Email : {\it contact@philippelefloch.org.}
\newline 
$^2$ CRIMERE, 14 rue du Cambodge, 75020 Paris, France. 
E-mail: {\sl jeanmarc.mercier@crimere.com}
}

\begin{abstract}
We revisit the method of characteristics for shock wave solutions to nonlinear hyperbolic problems 
and we describe a novel numerical algorithm ---the convex hull algorithm (CHA)--- in order to compute, both, 
entropy dissipative solutions (satisfying all relevant entropy inequalities) and entropy conservative (or multivalued) solutions to nonlinear hyperbolic conservation laws. Our method also applies to Hamilton-Jacobi equations and other problems endowed with a method of characteristics. From the multivalued solutions determined by the method of characteristic, our algorithm "extracts" the entropy dissipative solutions, even after the formation of shocks.
It applies to, both, convex or non-convex flux/Hamiltonians. We demonstrate the relevance of the proposed approach with a variety of numerical tests including a problem from fluid dynamics.  
\end{abstract}

\begin{keyword}  hyperbolic equations \sep
 method of characteristics \sep numerical algorithm  \sep compressible fluid

\

\MSC[2010] 76L05 \sep 65L05 \sep 35L60 \sep 	35L03 \hfill August 2014 
\end{keyword}

\end{frontmatter}



\section{Introduction}

The so-called `method of characteristics' allows one to determine {\sl smooth} solutions $u=u(t,x)$ to nonlinear hyperbolic equations. The solution formula takes the form $(u\circ S)(t,\cdot) = u_0 \circ S_0$ and is uniquely determined by the prescribed initial data of the Cauchy problem under consideration  
and, on the other hand, the family of characteristics $S=S(t,\cdot)$
associated with the problem.  (Here, $S_0 = S(0, \cdot)$ denotes a parametrization associated with the initial data.)  For nonlinear equations, this method generally fails for large times ---due to the formation of shock waves. 
In the present paper, we demonstrate that the characteristics are still of interest {\sl beyond the formation of shocks,} provided the method is suitably reformulated. 

\bei 

\item In our setting, it is natural to distinguish between two notions of solutions: 
\bei

\item  {\sl The entropy dissipative solutions,} which were defined in \cite{Hopf,Lax,Kruzkov} and satisfy all relevant entropy inequalities.

\item The {\sl entropy conservative solutions,} which satisfy {\sl entropy equalities} and are "highly oscillating" 
(and represent the multivalued solutions \cite{ForestierLeFloch}. 

\eei 

\item We exhibit here two novel formulas (see \eqref{u_b} and \eqref{u+}), which rely on a convex hull construction. Importantly, these formulas provide us with a 
numerical algorithm which is presented and implemented in the present paper. It is found to efficiently compute weak solutions with shocks for both convex and non-convex flux-functions. 

\eei 

We present our method in the context of multi-dimensional problems for, both, nonlinear hyperbolic 
conservation laws and  
Hamilton-Jacobi equations. We discuss several important properties of the solutions and perform numerical experiments. 
In particular, we demonstrate numerically that we recover the entropy solutions. The present paper is concerned with the numerical aspects of our approach, while a follow-up work will provide theoretical support for our new formulation. 
We also emphasize that the proposed method in principle should apply to a wide variety of nonlinear 
problems, provided they admit a method of characteristics. We apply it here to convex and non-convex 
conservation laws and to a problem of fluid dynamics. We expect our method to be applicable to other hyperbolic problems.  

We work in a spatial domain $\Omega$ which is taken to be $\Omega := \RR^D$ throughout this paper. We consider functions $u=u(t,x)$ defined for positive times $t\ge 0$
and satisfying one of the following two Cauchy problems: 
\begin{itemize}

\item {\it Nonlinear conservation laws:} 
\be \label{CL}
	\del_t u + \nabla \cdot f(u) = 0, \quad u(0,\cdot) = u_0,
\ee
where $f(u) := (f_d(u))_{d=1, \ldots, D} : \RR \mapsto \RR^D$ is a given flux-function 
and $\nabla \cdot f(u)$ denotes the divergence $\sum_{d=1, \ldots, D} \del_{x_d} f_d(u)$.

\item {\it Hamilton-Jacobi equations:}
$$
\del_t w +  H(\nabla w) = 0,
$$
with a given Hamiltonian $H : \RR^D \mapsto \RR$, where $\nabla u = \{ \del_{x_d} u \}_{d=1, \ldots, D} \in \RR^D$ denotes the gradient. 
In fact, by introducing the unknown 
$u := \nabla w \in \RR^D$, we can equivalently consider the following {\sl system} of nonlinear hyperbolic equations: 
\be \label{HJ}
	\del_t u +  \nabla \big( H(u) \big) = 0, \quad u(0,\cdot) = u_0 :=  \nabla w_0.
\ee
\end{itemize}
Clearly, the constraint $u = \nabla w$, i.e. the fact that $u$ is a gradient, is satisfied for all times provided it holds at the initial time. 

The method of characteristics applies to both problems and generates a smooth solution $u=u(t,\cdot)$ for 
sufficiently small times $t \le \delta$ (at least) in the following form (cf. the notation below): 
\be
u(t,\cdot)  = \big(u_0\circ S(0,\cdot) \big) \circ S^{-1}(t,\cdot), \quad  0\le t \le \delta,
\ee
where the map $S= S(t,x)$ parametrizes the family of characteristics. 
Our main observation based on the novel notion of {\sl convex hull transform}
allows us to generalize the characteristic method to all times, 
by relying on two explicit formulas (cf.~\eqref{u_b} and \eqref{u+}). 
Since the characteristic map $S(t, \cdot)$ is onto, defining its inverse $S^{-1}(t,\cdot)$ requires some care, 
and this issue is addressed below.  

Since our formulas are based on a reformulation of the characteristic method, the
approach we propose here is very general and should apply to many nonlinear hyperbolic problems. 
At the end of this paper, we for instance discuss an application to fluid dynamics. 
Our formulas are easily implemented in numerical applications (cf. the discussion in Section \ref{NR} below) and allows us to compute solutions for all times $t\ge 0$ and for all dimensions $D$. 
The method applies to arbitrary equations like \eqref{CL} and \eqref{HJ}, including \textsl{non-convex flux-functions} $f$ and \textsl{non-convex Hamiltonians}.  

In term of efficiency, our algorithm can favorably compete with spectral techniques such as the ones proposed for nonlinear hyperbolic problems by Sidilkover and Karniadakis \cite{SK} and followers, 
as well as techniques where the Rankine-Hugoniot relations are directly enforced at the discrete level such as the one proposed by Jaisankar and Raghurama Rao \cite{JR}. 
 Importantly, we expect the proposed method to generalize to many problems which admits a classical method of characteristics; for instance, problems with source-terms. Our approach could provide a robust alternative (or supplementary tool) to design well-balanced techniques, such the ones in Boscarino and Russo \cite{BoscarinoRusso}.


\section{The convex hull transformation} 

\subsection{Motivation} 

The present subsection provides a motivation for the present paper, but the reader unfamiliar with measure theory {\it 
can  skip this subsection} and read directly the next subsection. 
Without genuine loss of generality, we will work with smooth and spatially decaying initial data (so that the derivatives are integrable).   
Let $\Lambda \subset \RR^D$ be a convex set with unit Lebesgue measure and, typically, 
we can work with the unit cube $[-\frac{1}{2},\frac{1}{2}]^D$. 
By unit measure, we mean that $m(\Lambda)=1$, where $m$ stands for the Lebesgue measure. More generally, 
we use the notation $\mu(\Gamma)$ for the measure of a set $\Gamma$ with respect to a finite measure $\mu$.

Let us first consider general maps $S: y \in \Lambda \mapsto x=S(y) \in \Omega$, defined on an open, bounded, and convex set $\Lambda \subset \RR^D$ (with $D \geq 1$) and taking values in a convex and open set $\Omega \subset \RR ^D$. As usual, we denote by $\Pcal(\Omega)$ the set of all positive measures on $\Omega$ with unit total mass. 
Given two probability measures $\mu \in \Pcal(\Omega)$ and $\nu \in \Pcal(\Lambda)$, we say that
a map $S:\Lambda \mapsto \Omega$ is {\sl $(\mu, \nu)$--measure preserving} (or 
transports $\nu$ into $\mu$) if and only if 
$S_\sharp \nu = \mu$, that is, 
$\mu(S(A)) = \nu(A)$
for all ($\mu$-measurable) sets $A \subset \Omega$. 
This condition is equivalent to the change of variable formula
$\int_\Omega\varphi \, \mu = \int_\Lambda  (\varphi \circ S) \, \nu$.  

According to the general theory of optimal transport, given $\Omega \subset \RR^D$ and 
 for a positive probability measure $\mu \in \Pcal(\Omega)$ satisfying $\supp \mu = \Omega$ and $\mu << m$ (that is, $\mu$ is absolutely continuous 
with respect to the Lebesgue measure), these exists a unique map $S:\Lambda \mapsto \Omega$
''transporting'' the Lebesgue measure $m$ on $\Lambda$ into $\mu$.  
 More precisely, let $S:\Lambda \to \Omega$ be a square-integrable map  
and consider the probability measure $\mu := S_\# m \in \Pcal(\Omega)$. Then the following decomposition 
\be \label{PF}
	S= (\nabla h)\circ T, \qquad h: \Lambda \to \RR 
\text{ convex}, \qquad T_\#m = m,
\ee
holds and is called the \textsl{polar factorization} of the map $S$.  Here, by definition, the map $T$ is thus Lebesgue-measure preserving.

The following remark in also in order: the convex factorization is known to be unique when the measure 
$S_\#m := \mu > 0$ is absolutely continuous with respect to the Lebesgue measure, 
but uniqueness might be lost otherwise. In the present paper, due to the presence of shock waves, this uniqueness will not necessary hold. 
 

\subsection{A novel notion} 

Let $S:\Lambda \to \Omega$ be a square-integrable and surjective map
of the form $S = \nabla h$.
Then, the  {\it convex hull transformation} $S^+$ of $S$ is defined as
\be
	S^+:= \nabla h^+, \quad \text{where }h^+ \text{ is the convex hull of } h.
\ee
In view of this definition, $S^+$ is also surjective and the probability measure $\mu^+ := S^+_\# m \in \Pcal(\Omega)$ is in general singular, i.e. need not absolutely continuous with respect to the Lebesgue measure $m$ on $\Omega$ and may typically contain Dirac masses. In this context, we have the following proposition which allows us to inverse the map $S^+$ by 
solving the equation $u\circ S^+ = v_0$ when $v_0$ is a prescribed smooth function fixed throughout. 

\subsubsection*{Solving the equation $u\circ R = v_0$} 

Monotone maps admit generalized inverses defined as follows. 
Let $R:\Lambda \mapsto \Omega$ be a surjective map and consider the measure $\mu:=R_\#m>0$. 
Suppose $R=\nabla k$,  where $k: \Lambda \to \Omega$ is convex. 
Given an integrable function $v_0$ 
defined on $\Lambda$, we can introduce the integrable function
 $u(\cdot):=v_0 \circ R^{-1}(\cdot)$ 
 as the density of the measure $R_\# (v_0 m)$ with respect to $\mu$. This function is defined almost everywhere
 and satisfies for every set 
$A$ 
\be
	(u \mu)(A) = \int_B v_0, \quad B = R^{-1}(A). 
\ee
We omit here the discussion of the relevant regularity and functional spaces.  
With the convex hull transformation $S^+$,  the function $u = v_0 \circ (S^+)^{-1}$ defined as above, satisfies (by Jenssen inequality) 
\be
U( u ) \le U( v_0) \circ (S^+)^{-1} \qquad \text{ for every convex function } U
\ee
and, by definition, we then say that the convex hull transformation is {\it entropy dissipative.}  


\section{Entropy dissipative/conservative solutions} 

\subsection{The characteristic map} 
\label{schar}

\subsubsection*{Parametrization of the initial data}

We begin by formulating the standard characteristic method, which can be used to compute solutions for sufficiently 
small times. This holds for both equations \eqref{CL} and \eqref{HJ}. Observe that our formulation is based 
on a specific parametrization of the initial data. 

We consider either scalar-valued functions $u=u(t,x)$ 
or vector-valued maps $u=u(t,x) := \{u_i(t,x)\}_i$. As stated in the introduction, the gradient operator 
is denoted by $\nabla$, while the divergence operator is denoted by  $\nabla \cdot$). 
We introduce the norm $|\nabla u|_1:= \sum_i |\del_i u|$ for any (piecewise smooth, say) function
 $u=u(t, \cdot)$ defined on $\Omega$.  
If $u$ is vector-valued, we write $|\nabla u|_1:= \sum_{i,j} |\del_i u_j|$. 
With this notation,   
we are given a smooth initial data $u_0$ and we suppose $\supp \nabla u_0 = \Omega$ so that $u_0$ does not have ``flat part''. This is not a genuine restriction and can be ensured by an arbitrary perturbation of $u_0$. 
We then introduce the {\it initial map} 
\be \label{IM} 
\aligned
&
	S_0 = \nabla h_0 :\Lambda \mapsto \Omega, 
\qquad h_0: \Lambda \to \RR,  
\text{ convex},
\quad 
\\
& S_{0\#}m = \mu_0 = \frac{|\nabla u_0|_1}{|\nabla u_0|_1(\Omega)} \in \Pcal(\Omega).
\endaligned
\ee
We observe that $S_0$ is more regular (by one additional derivative) than the function $|\nabla u_0|_1$, which itself is
solely continuous 
(with 
bounded derivatives). 

\subsubsection*{Nonlinear conservation laws} 

For each time $t$, we define the {\it characteristic map} $S(t,\cdot):\Lambda \mapsto \Omega$ 
for the conservation law \eqref{CL} as :
\be \label{rhoeqref1}
    S(t,\cdot) := \big( S_0 + t f'( v_0 ) \big) 
\qquad v_0 := u_0 \circ S_0, 
\ee
which admits bounded derivatives up to first order at least. 

\subsubsection*{Hamilton-Jacobi equation}

Let $u_0: \Omega \to \RR^D$ 
be a smooth function, and suppose $\supp \nabla u_0 = \Omega$. 
Consider its initial map as in \eqref{IM}. Then define the {\it characteristic map} $S(t,\cdot):\Lambda \mapsto \Omega$ of the Hamilton-Jacobi equation \eqref{HJ} as :
\be \label{rhoeqref2}
    S(t,\cdot) = \big( S_0 + t (DH)( v_0 ) \big) 
\qquad v_0 := u_0 \circ S_0, 
\ee
where $DH$ is the differential of the Hamiltonian. 

Within a sufficiently small time interval $0 \le t \le \delta$, the characteristic maps above are invertible map, and provide
us with solutions to both equations, via 
\be \label{chars}
	u(t,\cdot) :=\big( v_0 \circ S^{-1}\big) (t,\cdot), \quad 0 \le t \le \delta.
\ee
Indeed, \eqref{chars} defines the unique smooth solution to \eqref{CL} and \eqref{HJ} within the interval under consideration. 


\subsection{A novel method of characteristics}

We now introduce two explicit formulas, which provide solutions to \eqref{CL} and to \eqref{HJ} for all positive times, i.e. beyond shock formation. 
Consider the characteristic maps $S(t,\cdot)$ for the conservation law and Hamilton-Jacobi equations defined in \eqref{rhoeqref1}-\eqref{rhoeqref2}, suppose without restriction that 
$S_0:\Lambda \to \Omega$ is square-integrable.  

\subsubsection*{Entropy conservative} 

Consider $S(t,\cdot) = \left(\nabla h \circ T \right) (t,\cdot)$ the polar factorization of $S$, see \eqref{PF}, with $h(t,\cdot): \RR^D \to \RR$ 
convex and $T : \Lambda \mapsto \Lambda$ Lebesgue measure preserving. For all positive times, we define the {\it entropy conservative} 
solution as
\be \label{u_b}
	u(t,\cdot) = \big(  v_0  \circ  T^{-1} \circ (\nabla h)^{-1} \big)(t,\cdot) 
\ee
for all positive times, in which the inverse $S^{-1} := T^{-1} \circ (\nabla h)^{-1}$ is well-defined since $h$ is convex and $T$ is measure-preserving.  

\subsubsection*{Entropy dissipative}

Suppose  $S(t,\cdot) = \nabla h(t,\cdot)$, with $h: \Omega \to \RR$ 
and consider its convex hull transform $S^+(t,\cdot)$, given earlier. 
We define  the {\it entropy dissipative solution} 
\be \label{u+}
	u(t,\cdot) = \big(  v_0 \circ  (S^+)^{-1}\big)(t,\cdot) 
\ee
for all positive times, in which the generalized inverse $S^{-1}$ was defined earlier. 
Furthermore, in the general situation $S = \nabla h \circ T$, we can 
introduce a more general transformation of the form  $S^+ = \nabla h^+ \circ T^+$, with a (possibly non-trivial)  measure-preserving map $T^+$, but this is unnecessary for our purpose in the present paper. 
 

\subsection{Properties of entropy dissipative solutions}

We now state some properties of the entropy dissipative solutions. Consider the map $S^+(t,\cdot)$ and the function $u(t,\cdot) = v_0 \circ (S^+)^{-1}$ defined by the  formula \eqref{u+}.
\begin{enumerate}

\item $S^+(t,\cdot): \Lambda \to \Omega$ is continuous and surjective. For almost all $x \in \Omega$, there exists a unique $y\in \supp (\nabla S^+)$ such that 
\be \label{XS}
	x=S^+(t,y)= S(t,y).
\ee

\item The function $u :=  v_0 \circ (S^+)^{-1}$ is a weak solution to the conservation law \eqref{CL}.

\item At every point of discontinuity of $u$ (with limits denoted by $u_l, u_r$), the map $S^+$ satisfies the Rankine-Hugoniot relations
$\del_t S^+ = \frac{ [ f( u \circ S^+ ) ] }{ [ u \circ S^+ ] } = \frac{ f( u_r) - f( u_l) }{u_r -u_l}$, 
while, at points of continuity, we can write 
$\frac{ [ f( u \circ S^+ ) ] }{ [ u \circ S^+ ] } = f'(u\circ S^+) = f'(v_0)$ (in the support $\nabla S$). 

\item The entropy condition holds: for instance in one space dimension and at a point of discontinuity where
$u_l \leq u_r$ one has 
$\frac{f(u) - f(u_l)}{u-u_l} \ge \frac{f(u_r) - f(u_l)}{u_r-u_l}\ge  \frac{f(u) - f(u_r)}{u-u_r}$
for all $u\in[u_l,u_r]$.

\item In one space dimension and for a convex flux-function, we can recover the Hopf-Lax formula: 
\eqref{u+} coincides with the following formula $w(t,x) = \int_{-\infty}^x u(t,s) \, ds$ with 
$w(t,x) = \inf_z \big( w_0(z) + t f^*\big(\frac{x-z}{t}\big) \big)$, 
where $f^*(z) = \sup_{y}\left(zy-f(y)\right)$ is the Legendre-Fenchel transform of the flux.
\end{enumerate}

We recall that the 
Hopf-Lax formula allows one to compute solutions to (one-dimensional) 
conservation laws and (multi-dimensional) Hamilton-Jacobi equations when the 
flux or the Hamiltonian are convex.  


\section{An illustration with Burgers equation}  
\label{NR}

Before we present our algorithm in full details (in the forthcoming section), we want first to illustrate our formula with a typical example. 
We thus consider the formulas \eqref{u+}-\eqref{u_b} for the inviscid Burgers equation, that is, the equation
 \eqref{CL} with quadratic flux-function: 
\be
	\del_t u + \del_x \Big( {u^2 \over 2} \Big) = 0, \qquad u(0,x) = u_0(x).
\ee 
Here, we define the map $S : \Lambda=[-{1 \over 2}, {1 \over 2} ] \to \Omega := \RR$ by 
\be
S(t, \cdot) = S_0 + t v_0, \qquad S_{0,\#} m = \frac{ |\del_x u_0|}{TV(u_0)}, \qquad v_0 := u_0 \circ S_0.
\ee 

\begin{figure}[h]
\centering
\begin{tabular}{ccc}

\begin{minipage}{1.3in}
\includegraphics[height=1.6in,width=1.9in]
{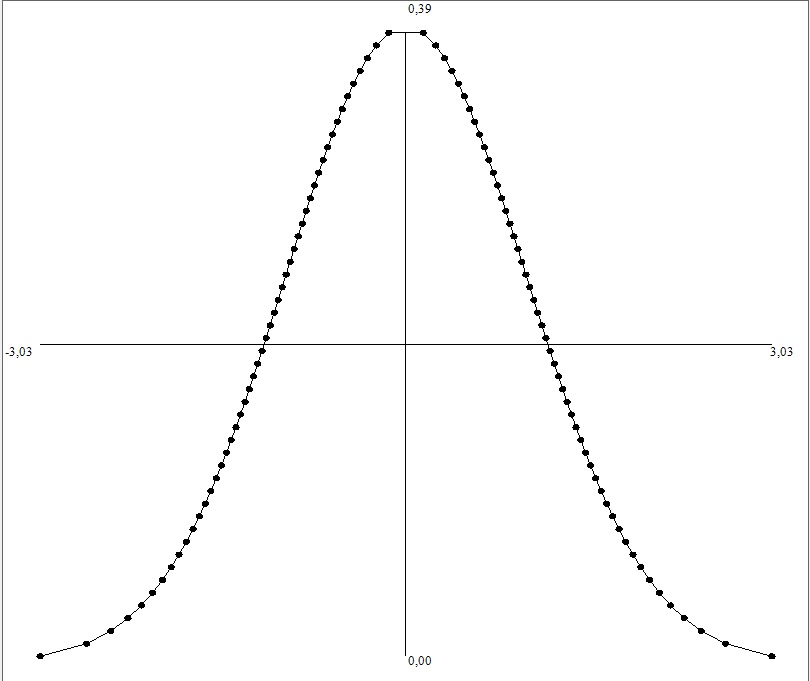}
{(a) t=0: initial condition}
\end{minipage}

\hskip3.01cm

\begin{minipage}{1.5in}
\includegraphics[height=1.6in,width=1.9in]
{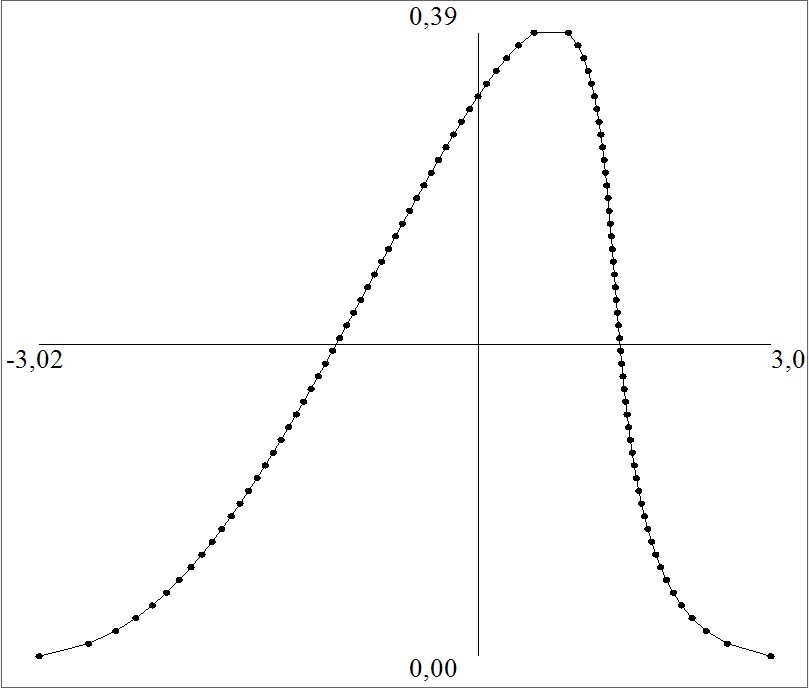}
(b) t=3: smooth regime
\end{minipage}

\\

\\

\begin{minipage}{1.5in}
\includegraphics[height=1.6in,width=1.9in]
{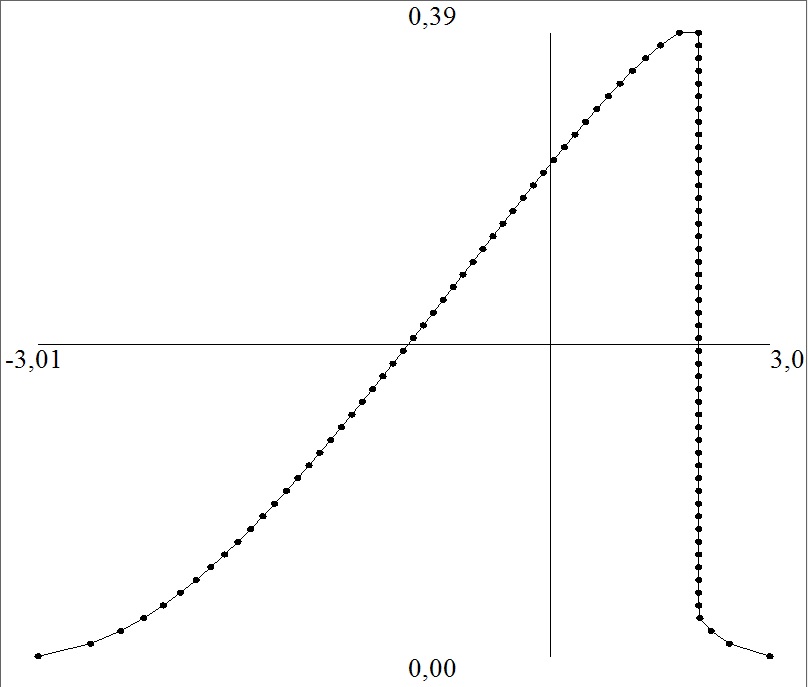}
(c)t=6: shock formation
\end{minipage}

\hskip3.01cm

\begin{minipage}{1.5in}
\includegraphics[height=1.6in,width=1.9in]
{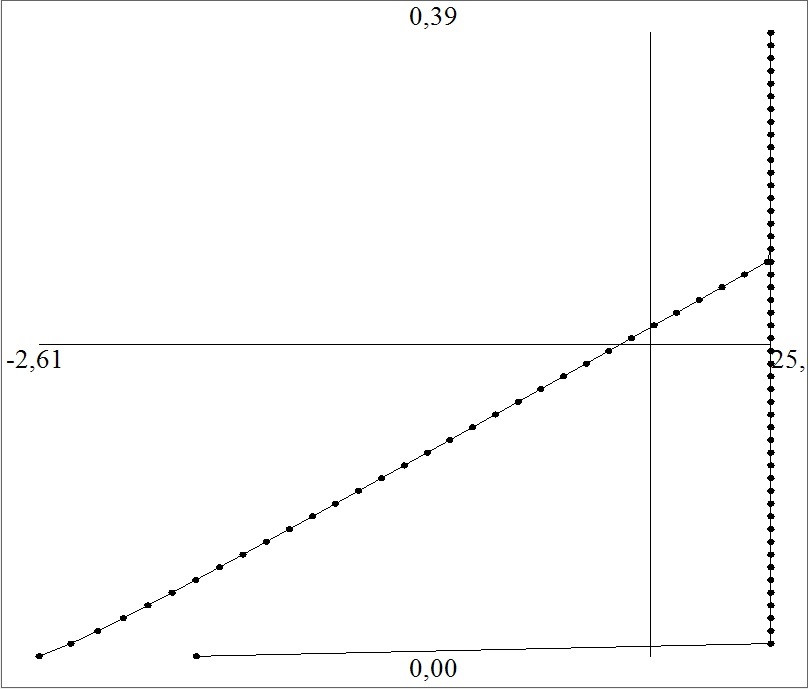}
(d) t=106: N-wave regime.
\end{minipage}
\end{tabular}
\caption{Entropy dissipative solutions to Burgers equation} \label{BEE}
\end{figure}

Our numerical results in Figure \ref{BEE} plots the solution computed using the formula \eqref{u+} 
with the Gaussian initial data
\be
	 u(0,x) = {1 \over { \sqrt{ 2 \pi }}} \exp\Big( - {x^2 \over 2 } \Big).
\ee
As will be explained in the following section, the solution $u=u(t,x)$ is represented by a cloud of points $\big( S^+(t,y_i), u\circ S^+(t,y_i) \big)_{i=0, \ldots, N-1}$, with $y_i = \frac{i+1/2}{N}-{1 \over 2}$. The simulation is ran with $N=100$. 

Figure \ref{BEE} (a) represents the initial Gaussian condition. Figure \ref{BEE} (b) shows the solution $u(t=3,\cdot)$ during the smooth regime of evolution. Next, Figure \ref{BEE} (c) 
displays the solution after the shock formation. 
Finally, Figure \ref{BEE} (d) (at time $t=106$) demontrates the convergence of the solution toward the so-called N-wave, as expected, 
which is the asymptotic profile solution $\frac{x}{\sqrt{t}}$ (up to a translation and to a jump). Observe the
 formation of a "spike" in the shock, which is clearly visible at the time $t=106$. 
Namely, our method generates an additional spike within the jump discontinuity (which could be easily filtered and removed, if one wishes). 
Note also the isolated point in Figure \ref{BEE} (d), which is a minor numerical artifact.

\begin{figure}[h]
\centering
\begin{tabular}{ccc}

\begin{minipage}{1.3in}
\includegraphics[height=1.6in,width=1.9in]
{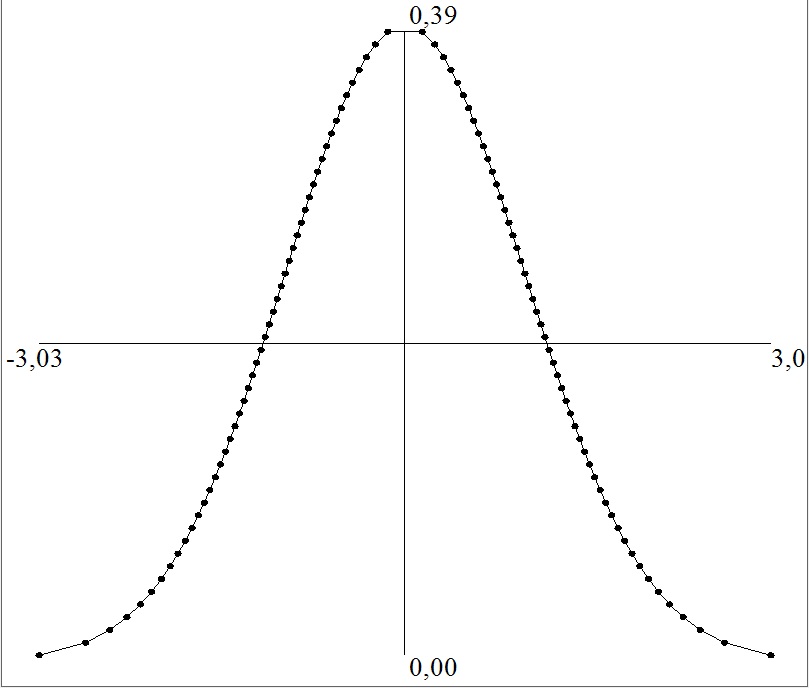}
(a) t=0: initial condition
\end{minipage}

\hskip3.01cm

\begin{minipage}{1.5in}
\includegraphics[height=1.6in,width=1.9in]
{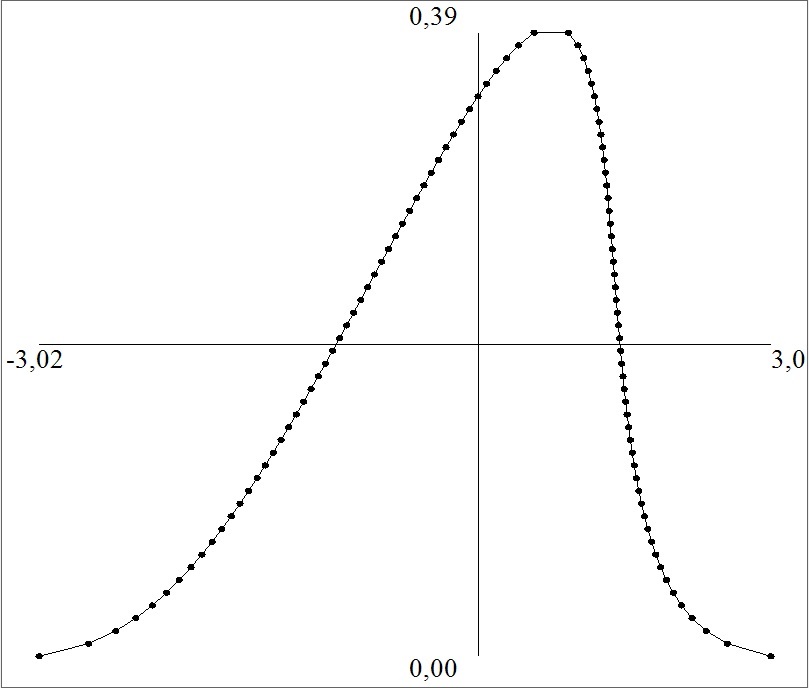}
(b) t=3: smooth regime 
\end{minipage}

\\

\begin{minipage}{1.5in}
\includegraphics[height=1.6in,width=1.9in]
{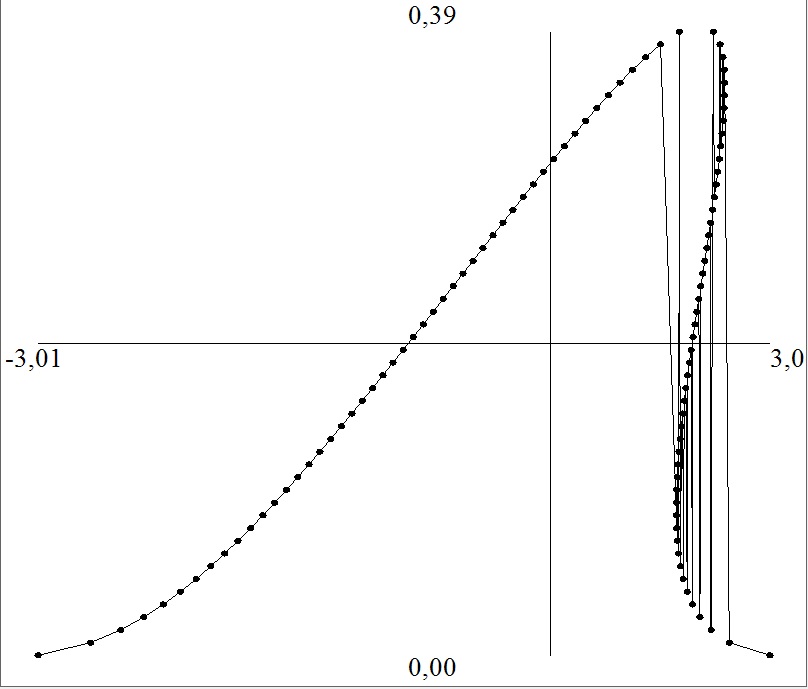}
(c) t=6: oscillating regime
\end{minipage}

\hskip3.01cm

\begin{minipage}{1.5in}
\includegraphics[height=1.6in,width=1.9in]
{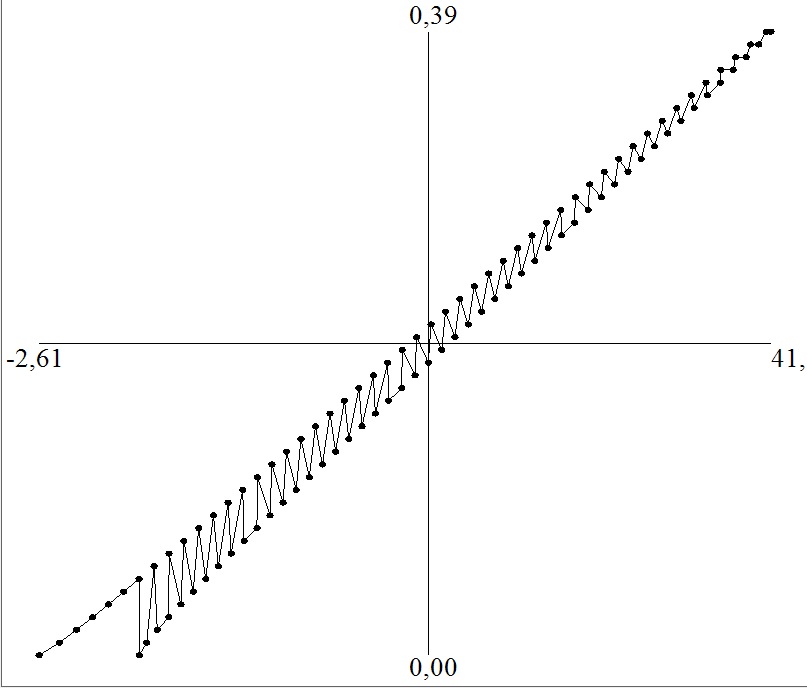}
(d) t=106: long time behavior
\end{minipage}
\end{tabular}
\caption{Entropy conservative solutions to Burgers equation} \label{BEC}
\end{figure}

In contrast, Figure \ref{BEC} represents the {\it conservative solution} with the same initial condition, but
 computed with our formula \eqref{u_b}. We can easily deduce this solution from the multivalued one. 
The solution is selected among all possible values by using the re-ordering map $T$ in the polar factorization \eqref{PF}.

\section{Implementation of the convex hull algorithm (CHA)}
\label{algo}

We now present our algorithm for one-dimensional conservation laws of the form \eqref{CL}, that is, 
\be
\label{eq501} 
	\del_t u + \del_x f(u) = 0, \qquad u(0,x) = u_0(x), \quad x \in \Omega \subset \RR.  
\ee
We rely on our formulas \eqref{u+} and \eqref{u_b}, proposed in the first part of this paper.  
We introduce the map $S : \Lambda \mapsto \Omega := \RR$ by 
\be
S(t, \cdot) = S_0 + t f'(v_0), \quad S_{0,\#} m = \frac{ |\del_x u_0|}{TV(u_0)}, \quad v_0 := u_0 \circ S_0.
\ee

We are going to our numerical algorithm, step by step, by focusing on Burgers equation, for simplicity in the presentation. The proposed algorithm yields the solution \eqref{u+} or \eqref{chars}, and works equally well in any spatial dimension and and with Hamilton-Jacobi equations.  We point out that the code based on our method is  available \cite{Mercier}, so that our results can be easily reproduced by the reader. This algorithm will be tested with more complex case studies in the following section.

\begin{figure}[h]
\centering
\begin{tabular}{ccc}

\begin{minipage}{1.3in}
\includegraphics[height=1.6in,width=1.9in]
{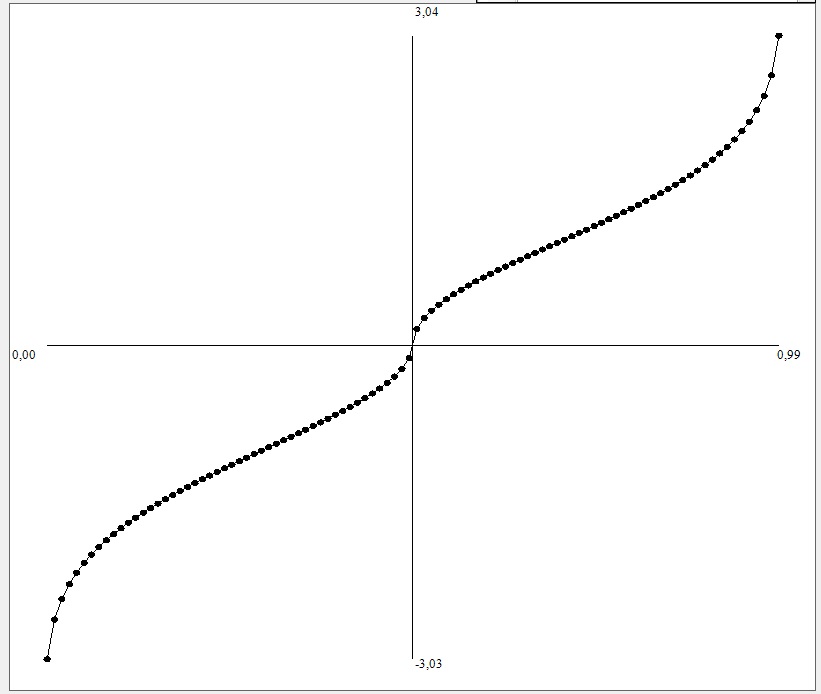}
{(a) t=0}
\end{minipage}

\hskip3.01cm

\begin{minipage}{1.3in}
\includegraphics[height=1.6in,width=1.9in]
{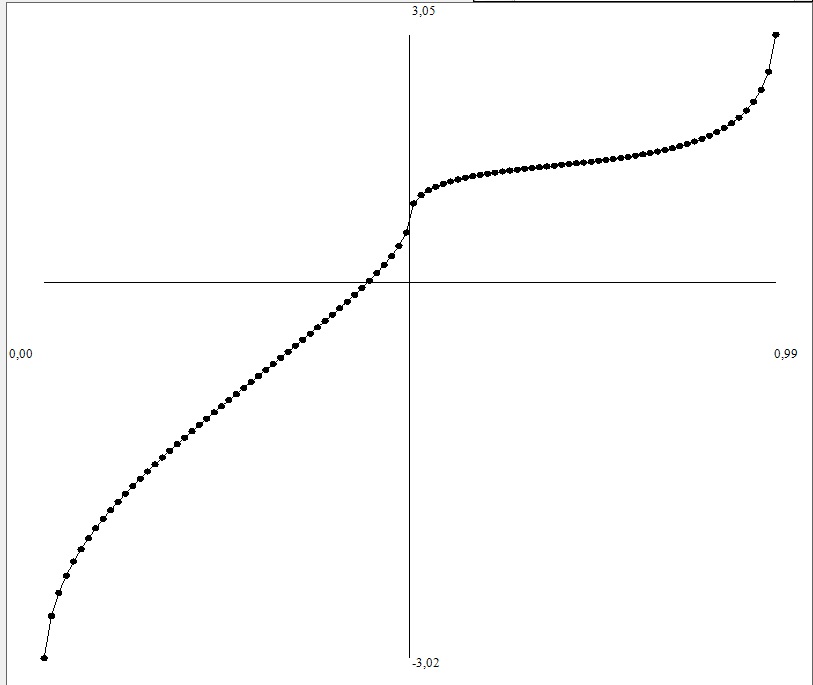}
{(b) t=3}
\end{minipage} 

\\

\begin{minipage}{1.3in}
\includegraphics[height=1.6in,width=1.9in]
{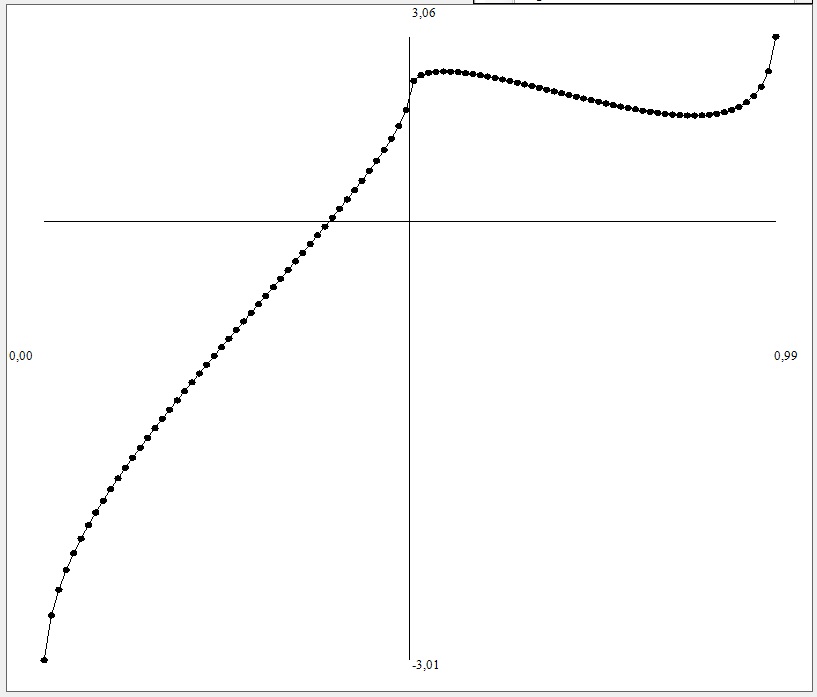}
{(c) t=6}
\end{minipage}

\hskip3.01cm

\begin{minipage}{1.3in}
\includegraphics[height=1.6in,width=1.9in]
{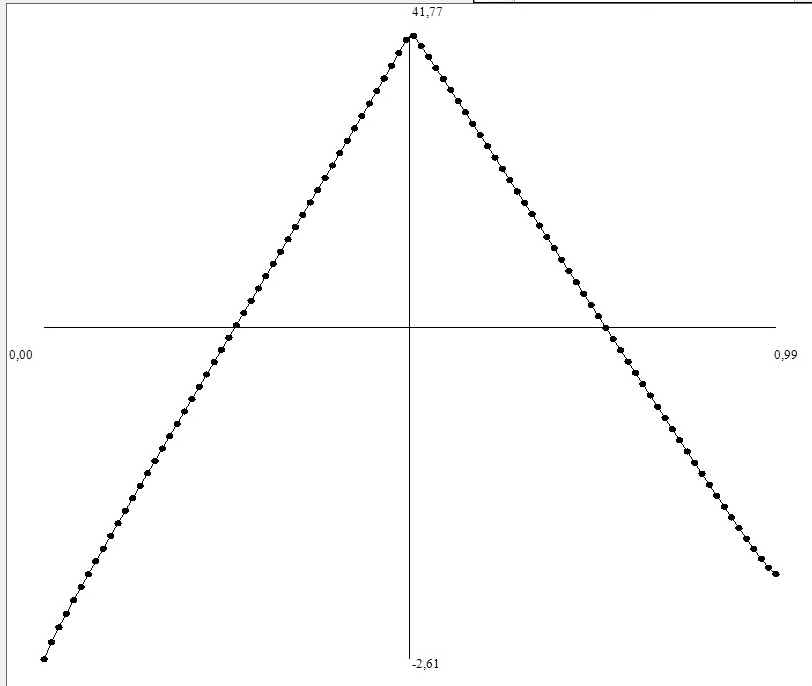}
{(d) t=106}
\end{minipage}

\end{tabular}
\caption{The characteristic map $S(t,\cdot)$} \label{SE}
\end{figure}

\begin{figure}[h]
\centering
\begin{tabular}{ccc}

\begin{minipage}{1.3in}
\includegraphics[height=1.6in,width=1.9in]
{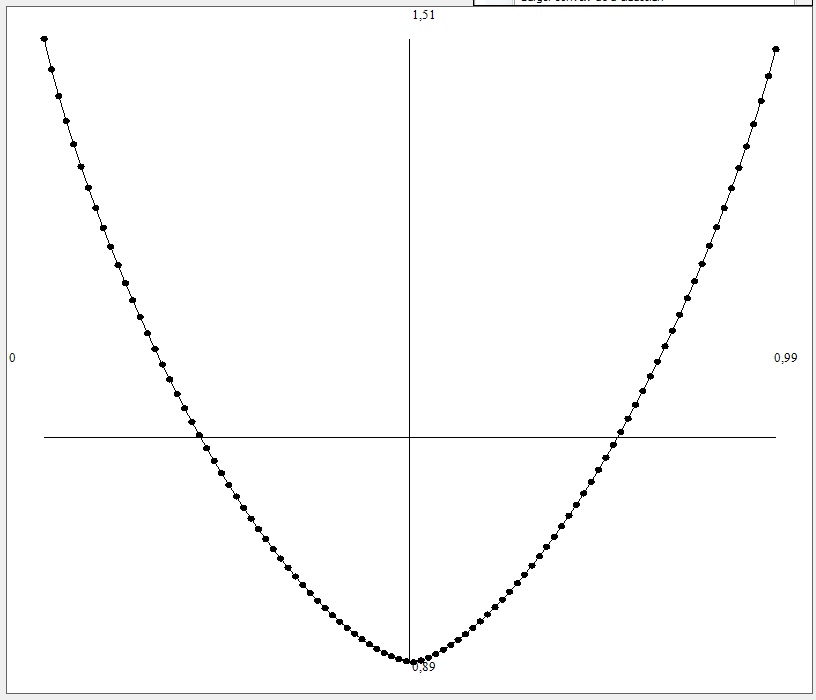}
(a) t=0
\end{minipage}

\hskip3.01cm

\begin{minipage}{1.3in}
\includegraphics[height=1.6in,width=1.9in]
{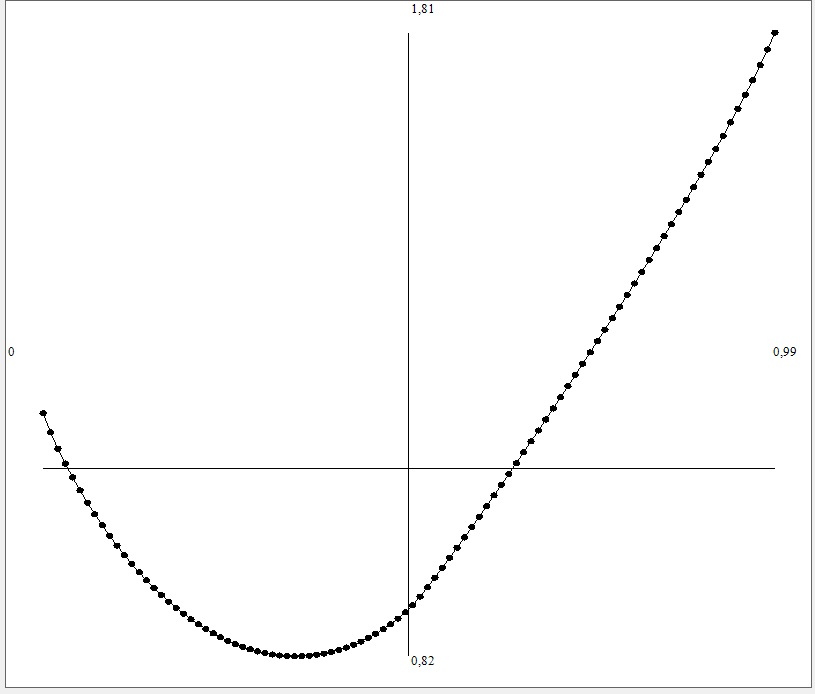}
(b) t=3
\end{minipage}

\\

\begin{minipage}{1.3in}
\includegraphics[height=1.6in,width=1.9in]
{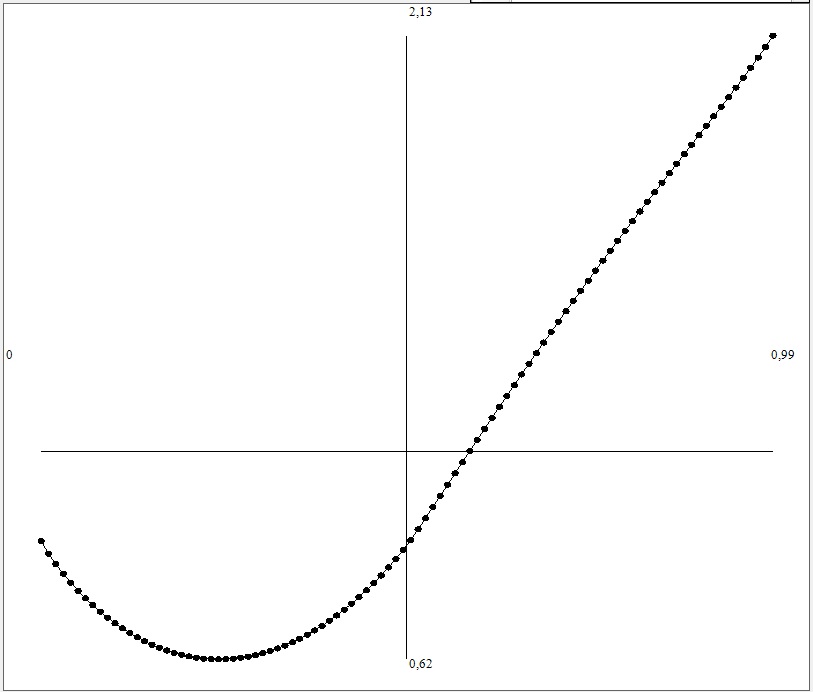}
(c) t=6
\end{minipage}

\hskip3.01cm

\begin{minipage}{1.3in}
\includegraphics[height=1.6in,width=1.9in]
{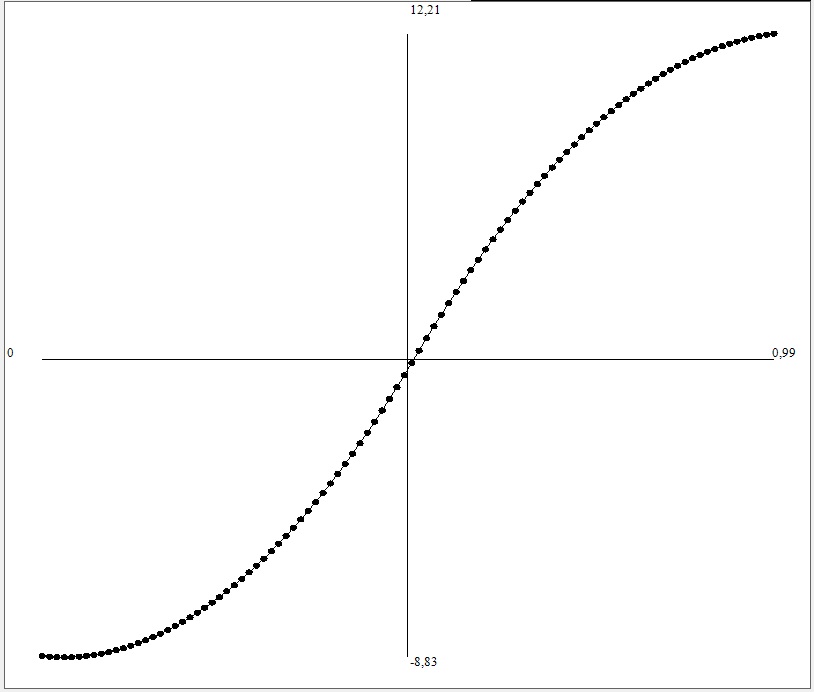}
(d) t=106
\end{minipage}

\end{tabular}
\caption{The map $h(t,\cdot)$} \label{hE}
\end{figure}

\begin{figure}[h]
\centering
\begin{tabular}{ccc}

\begin{minipage}{1.3in}
\includegraphics[height=1.6in,width=1.9in]
{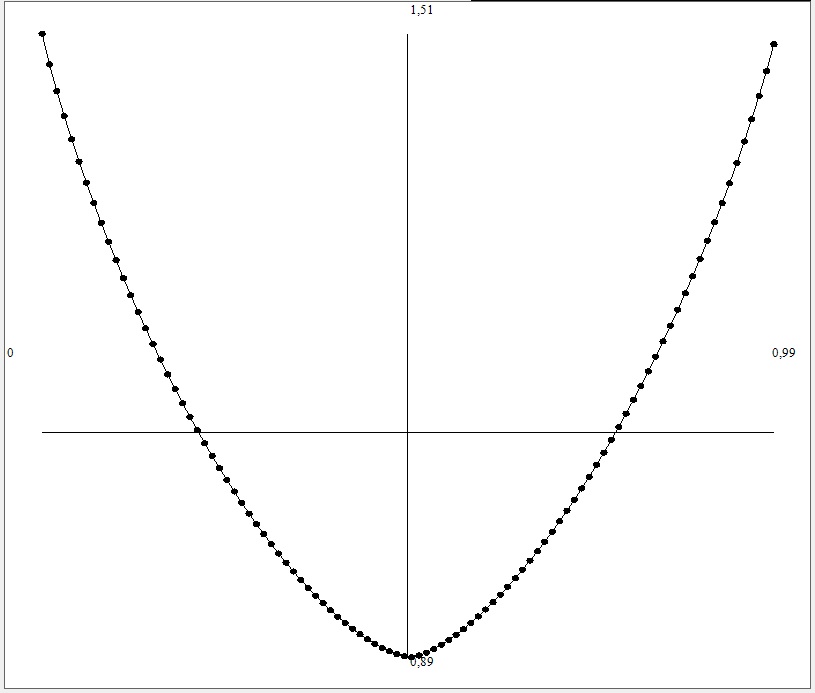}
(a) t=0
\end{minipage}

\hskip3.01cm

\begin{minipage}{1.3in}
\includegraphics[height=1.6in,width=1.9in]
{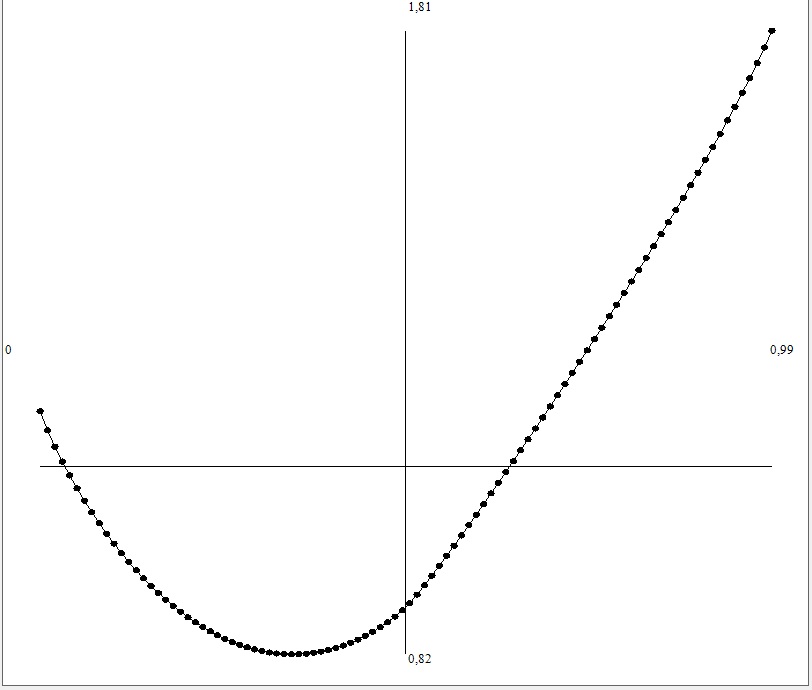}
(b) t=3
\end{minipage}

\\

\begin{minipage}{1.3in}
\includegraphics[height=1.6in,width=1.9in]
{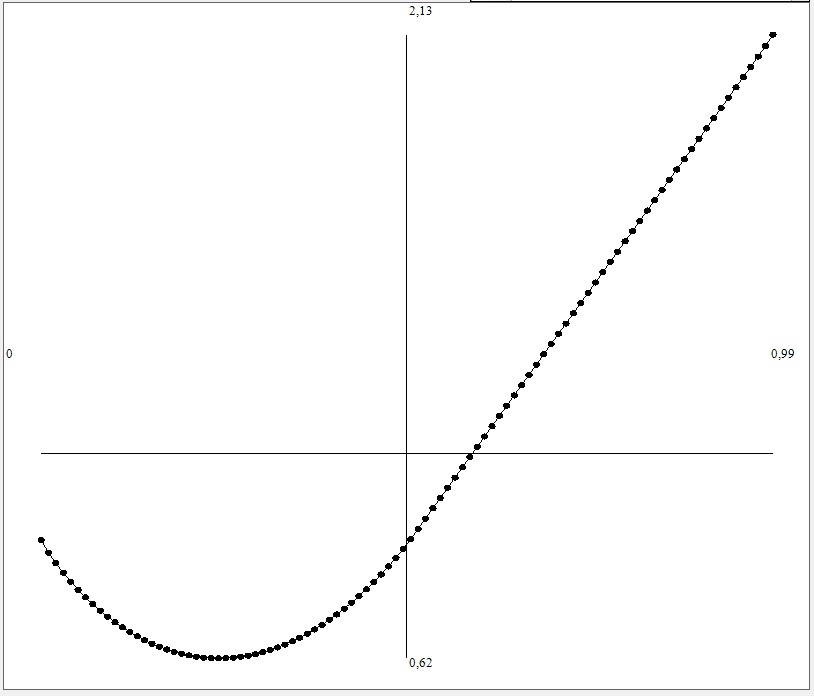}
(c) t=6
\end{minipage}

\hskip3.01cm

\begin{minipage}{1.3in}
\includegraphics[height=1.6in,width=1.9in]
{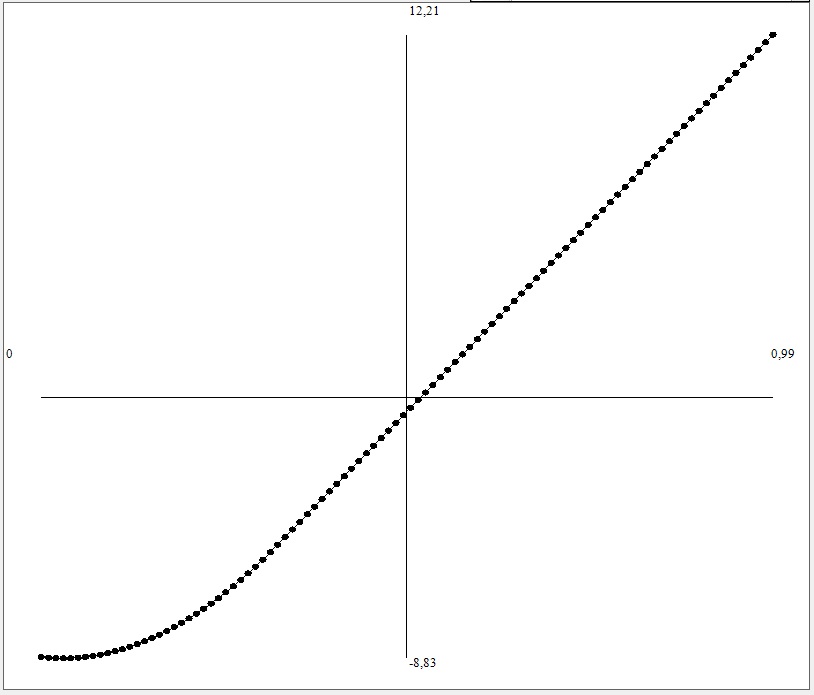}
(d) t=106
\end{minipage}

\end{tabular}
\caption{Convex hull transform component $h^+(t,\cdot)$} \label{hE+}
\end{figure}

\begin{figure}[h]
\centering
\begin{tabular}{ccc}

\begin{minipage}{1.3in}
\includegraphics[height=1.6in,width=1.9in]
{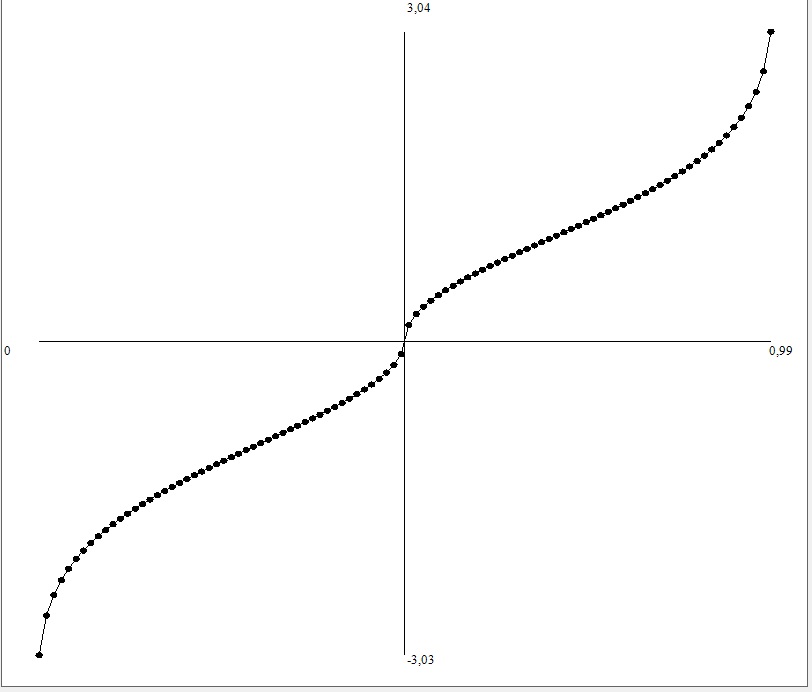}
{(a) t=0}
\end{minipage}

\hskip3.01cm

\begin{minipage}{1.3in}
\includegraphics[height=1.6in,width=1.9in]
{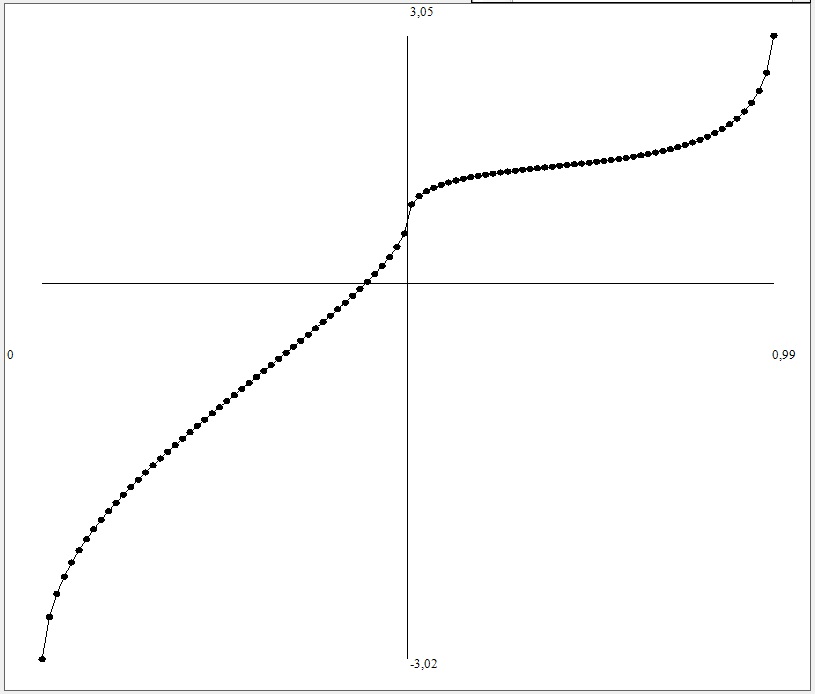}
{(b) t=3}
\end{minipage}

\\

\begin{minipage}{1.3in}
\includegraphics[height=1.6in,width=1.9in]
{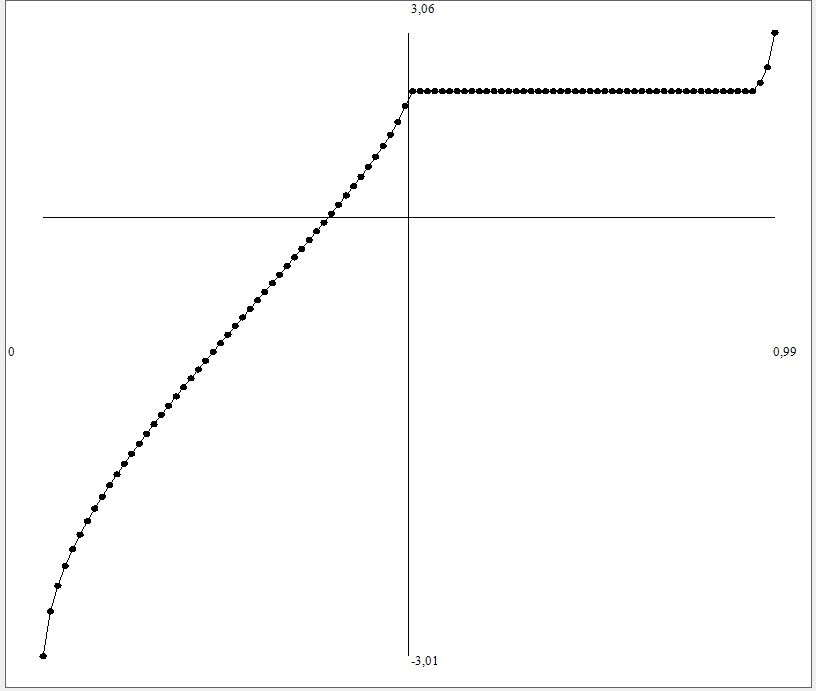}
{(c) t=6}
\end{minipage}

\hskip3.01cm

\begin{minipage}{1.3in}
\includegraphics[height=1.6in,width=1.9in]
{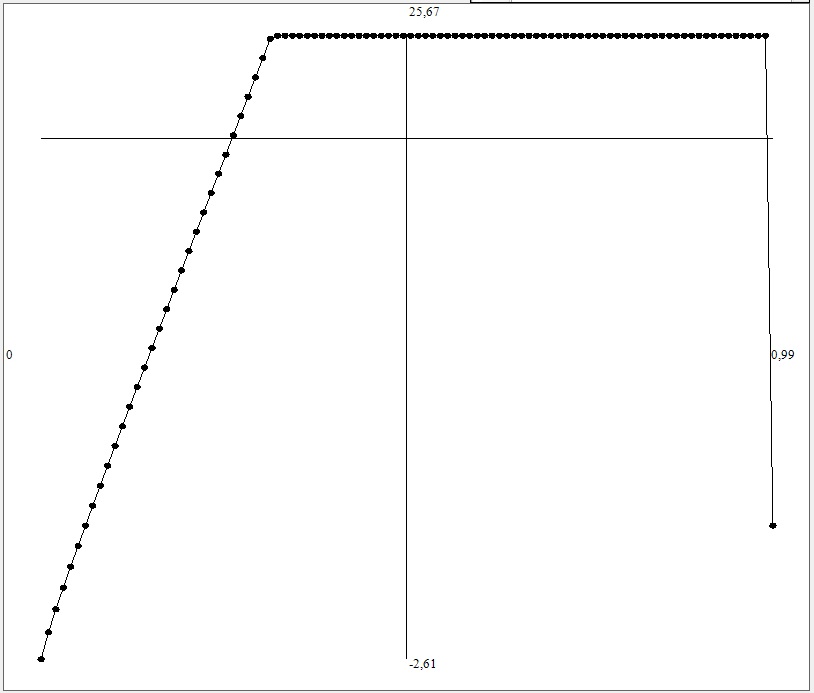}
{(d) t=106}
\end{minipage}

\end{tabular}
\caption{Convex hull transform $S^+(t,\cdot)$} \label{SE+}
\end{figure}

We focus on the formula \eqref{u+} and on Burgers equation studied in the previous section. Our construction is as follows: 
we consider the characteristic map $S(t,\cdot) = S_0+tf'(v_0) : \Lambda \mapsto \Omega$, where $v_0:=u_0\circ S_0$
and we set $S(t,\cdot) := \nabla h(t,\cdot)$. For the numerical tests, we will always use $\Lambda = (-\frac{1}{2},\frac{1}{2})^D$ (with $D=1$ in this paper) and $\Omega \subset \RR^D$. 
Our {\it Convex Hull Algorithm} (CHA) amounts to determine the solution 
\be
	u(t,\cdot) = \big(  v_0 \circ  (S^+)^{-1}\big)(t,\cdot), 
\ee
in which 
\be
S^+(t,\cdot) = \nabla h^+(t,\cdot), \quad h^+(t,\cdot) \text{ convex hull of } h(t,\cdot). 
\ee
The main steps of this algorithm are as follows, for each time $t$ :
\begin{itemize}

\item {\it Step 1.} Compute the characteristic map $S(t,\cdot)$ (as illustrated by Figure \ref{SE}) and its companion function $h(t,\cdot)$ (see Figure \ref{hE}), deduced coming from of a Helmholtz-Hodge decomposition of the form 
\be
S(t,\cdot) = \big( \nabla h + \chi \big)(t,\cdot),
\ee
where the map $\chi$ is a divergence-free, i.e. $\nabla \cdot \chi = 0$. In one space dimension, we can always normalize $\chi$ to vanish identically. 

\item {\it Step 2.} Compute the convex hull $h^+(t,\cdot)$ of the function $h(t,\cdot)$. This step is illustrated by Figure \ref{hE+}.

\item {\it Step 3.} Compute again the map $S^+(t,\cdot) = \nabla h^+(t,\cdot)$ (see Figure \ref{SE+}) as well
as the composite function $\big( u \circ S^+ \big)(t,\cdot) := v_0$ (as represented in Figure \ref{BEE}).
\end{itemize}

Concerning Step 2, efficient algorithms in order to compute convex hulls are available as open source softwares.
(See, for instance, http://www.cgal.org.)  We used the specialized QHull library available at http://www.qhull.org. Our numerical contribution for the present work is a Helmoltz-Hodge decom\-position/re-composition code based
on multi-dimensional unstructured meshes, which can be described as follows.

Let us introduce the following non-local basis functions $\psi_Y(X) : \Lambda \mapsto \RR$, where $Y = (y_d )_{d=1, \ldots, D} \in \Lambda$ :
\be \label{psi}
	\psi_{Y}(X) = \frac{1}{2} |X-Y|_1 = \frac{1}{2} \sum_{d=1, \ldots, D} |x_d-y_d|,\quad X\in \Lambda.
\ee
Let ${\it Y} := \big(Y_n \big)_{n=1, \ldots, N}$, $N$ be distincts sampling points of the computational 
domain $\Lambda$, with $Y_j := \frac{i+1/2}{N}-{1 \over 2}$ in our one-dimensional setting. 
We then consider any (possibly) vector-valued map $S : \Lambda \mapsto \Omega$ and the following minimization problem
\be
	\inf_{(\alpha_n)_{n=1, \ldots, N} \in \RR^N} \int_\Lambda |S- \sum_{n = 1, \ldots, N} \alpha_n \nabla \psi_{Y_n}|^2_2 \, dm,
\ee
which can be interpreted as a classical projection method. This latter problem is approximated by a linear system corresponding to a Gram-Schmidt algorithm (in which the notation $\langle f, h \rangle_{L^2(\Lambda, \RR^D)}$ 
is used for the integral $\int_\Lambda \langle f_1, f_2 \rangle$ of the scalar product of 
two vector-valued functions):  
\be
\aligned
&	\sum_{j= 1, \ldots, N} \alpha_n  \langle \nabla \psi_{Y_i},\nabla  \psi_{Y_j} \rangle_{L^2(\Lambda,\RR^D)} 
  = \frac{1}{2} \langle S, \nabla \psi_{Y_i} \rangle_{L^2(\Lambda, \RR^D)} 
\\
&\simeq \frac{1}{2N} \sum_{n=1, \ldots, N \atop n \neq i} S(Y_n) \cdot \nabla \psi_{Y_i}(Y_n), \quad i =1, \ldots, N. 
\endaligned
\ee
This linear system is invertible, provided $Y_i \neq Y_j, i\neq j$, and the right-hand side approximation makes sense 
for all maps $S$ that are finite sums of convex and concave functions. Once the coefficients $\alpha_n$
 are computed, the component $h_{\it Y}$ is obtained by 
\be
h_{\it Y} = \sum_{n=1, \ldots, N} \alpha_n \psi_{Y_n}.
\ee
 Once the convex hull $h^+_{\it Y}$ is computed, a similar algorithm yields us the components 
$\alpha^+_n$ of the projection $h^+_{\it Y} = \sum_n \alpha^+_n \psi_{Y_n}$. This allows us to finally compute 
\be
S^+_{\it Y} = \sum_n \alpha^+_n \nabla \psi_{Y_n}.
\ee


\section{Applications}

\subsection{A fluid dynamical problem} 
 
Finally, we consider the coupling between Burgers equation satisfied by the velocity function $u$, that is, 
\be
\del_t u + \del_x \Big(\frac{u^2}{2} \Big) = 0
\ee
and a linear transport equation for a passive scalar $w$ (cf.~\cite{LeFloch90} for a mathematical background): 
\be
\del_t w + u \del_x w = 0. 
\ee
By setting $\rho := w_x$, we can interpret $(u,w)$ as solutions to the system of fluid dynamics  
\be \label{345}
\partial_t \rho+\partial_x (\rho u)=0, 
\ee
\be \label{346}
\partial_t (\rho u)+\partial_x (\rho u^2)=0. 
\ee  

Let $S=S(t,\cdot)$ be the characteristic map for Burgers equation, i.e. $S_t = v_0 = u \circ S$. Consider the function 
\be
w(t,\cdot) = w(0,\cdot) \circ S(0,\cdot) \circ S(t,\cdot)^{-1}
\ee
with which we can compute (at least formally) 
$
	0=\big( w\circ S \big)_t = w_t \circ S + S_t w_x \circ  S = \big( w_t +  w_x  u \big) \circ S.
$
Hence,  $w_t +  w_x  f'(u) = 0$, provided $S$ is surjective, so that \eqref{345}-\eqref{346} is satisfied. Consequently, we obtain the  function associated with the density $w(t,\cdot) = w(0,\cdot) \circ S(0,\cdot) \circ S^{-1}(t,\cdot)$.

\begin{figure}[h]
\centering
\begin{tabular}{ccc}

\begin{minipage}{1.3in}
\includegraphics[height=1.6in,width=1.9in]
{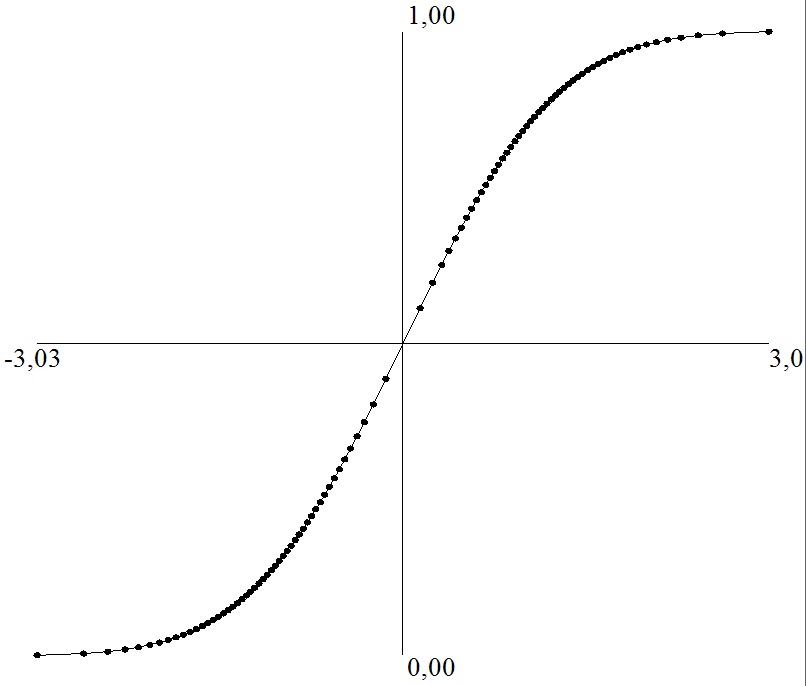}
{(a) t=0}
\end{minipage}

\hskip3.01cm

\begin{minipage}{1.3in}
\includegraphics[height=1.6in,width=1.9in]
{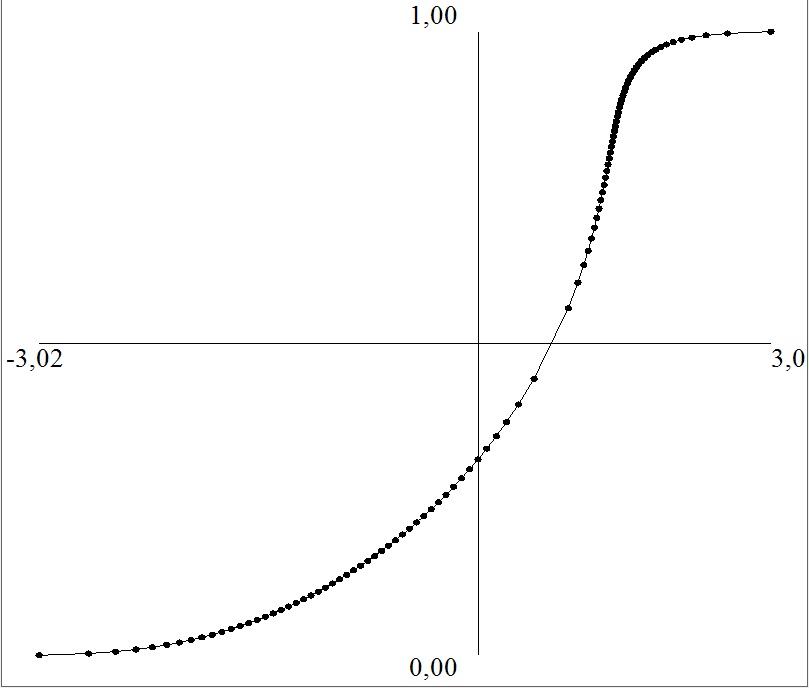}
{(b) t=3}
\end{minipage}

\\

\begin{minipage}{1.3in}
\includegraphics[height=1.6in,width=1.9in]
{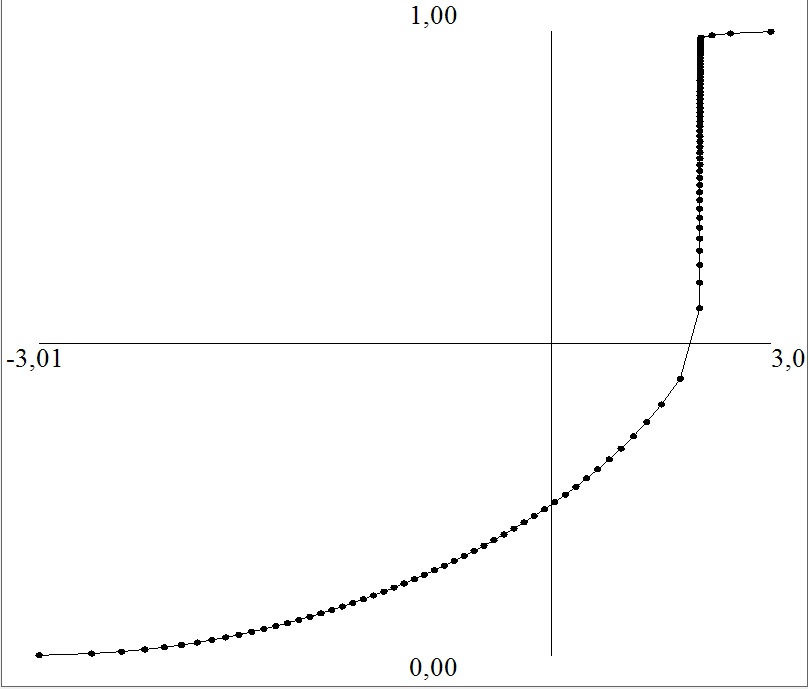}
{(c) t=6}
\end{minipage}

\hskip3.01cm

\begin{minipage}{1.3in}
\includegraphics[height=1.6in,width=1.9in]
{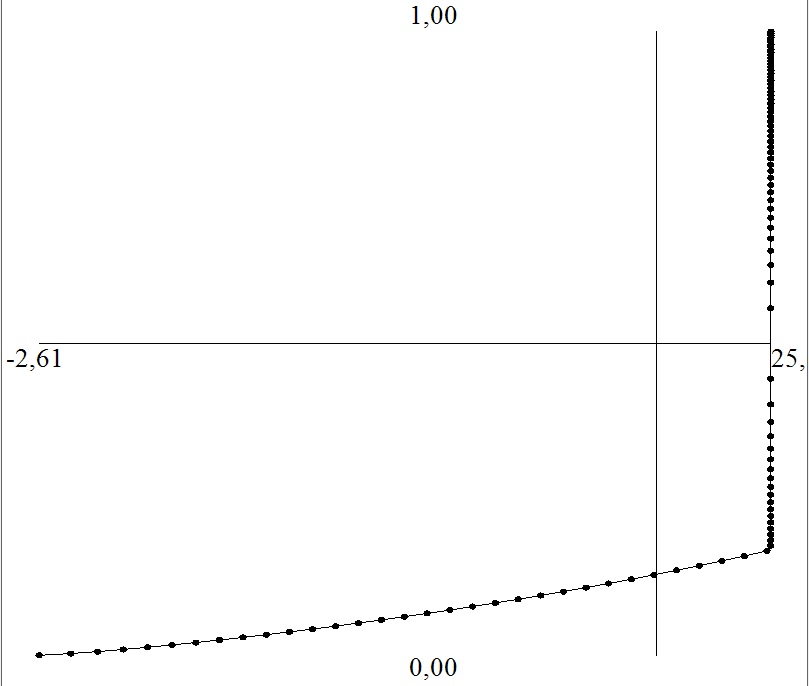}
{(d) t=106}
\end{minipage}

\end{tabular}
\caption{Solution component $w(t,\cdot)$} \label{Euleru}
\end{figure}

\begin{figure}[h]
\centering
\begin{tabular}{ccc}

\begin{minipage}{1.3in}
\includegraphics[height=1.6in,width=1.9in]
{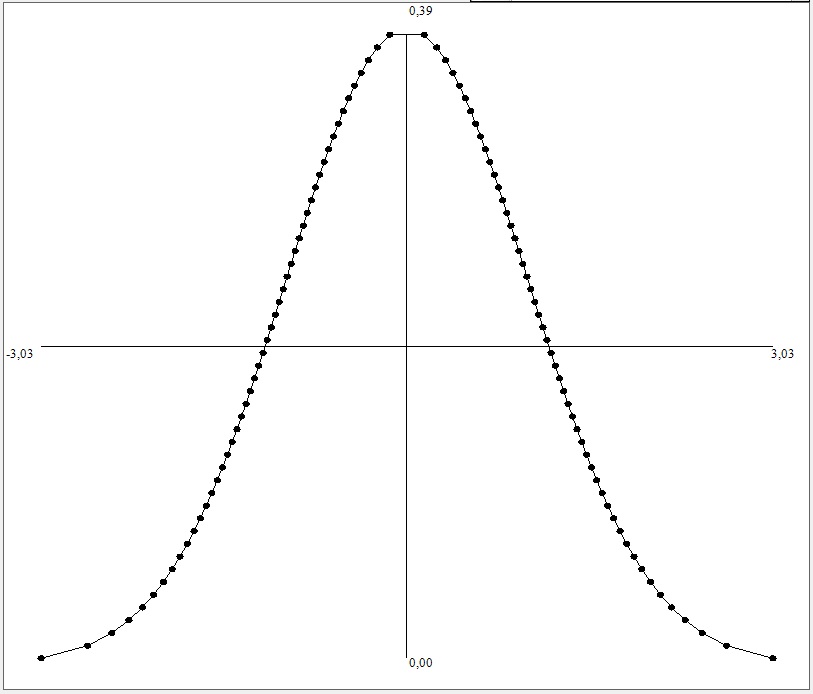}
{(a) t=0}
\end{minipage}

\hskip3.01cm

\begin{minipage}{1.3in}
\includegraphics[height=1.6in,width=1.9in]
{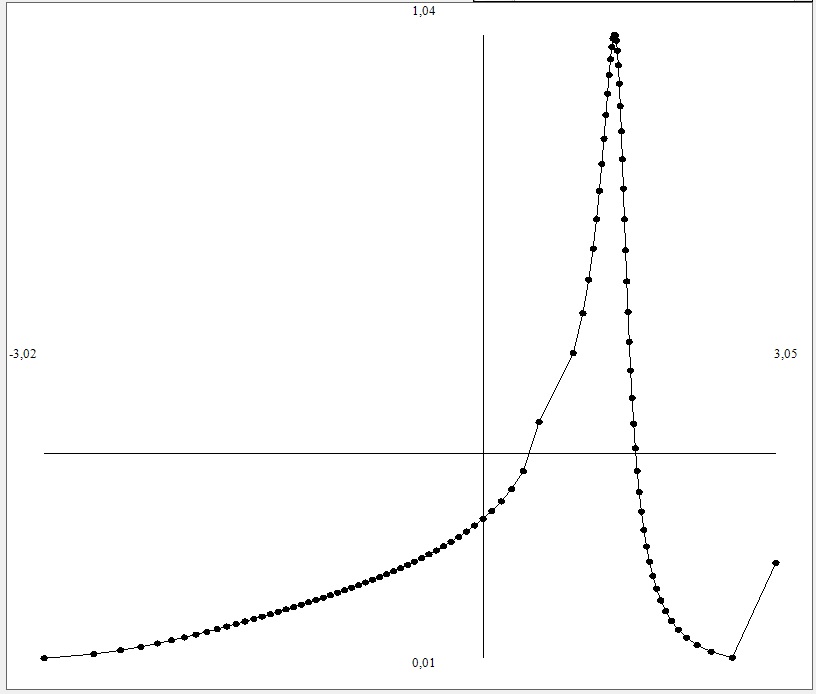}
{(b) t=3}
\end{minipage}

\\

\begin{minipage}{1.3in}
\includegraphics[height=1.6in,width=1.9in]
{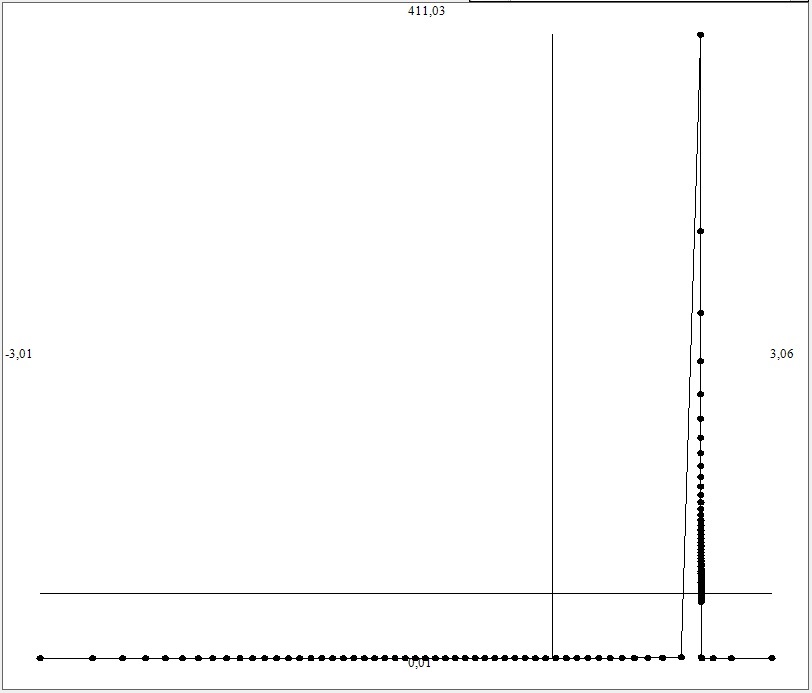}
{(c) t=6}
\end{minipage}

\hskip3.01cm

\begin{minipage}{1.3in}
\includegraphics[height=1.6in,width=1.9in]
{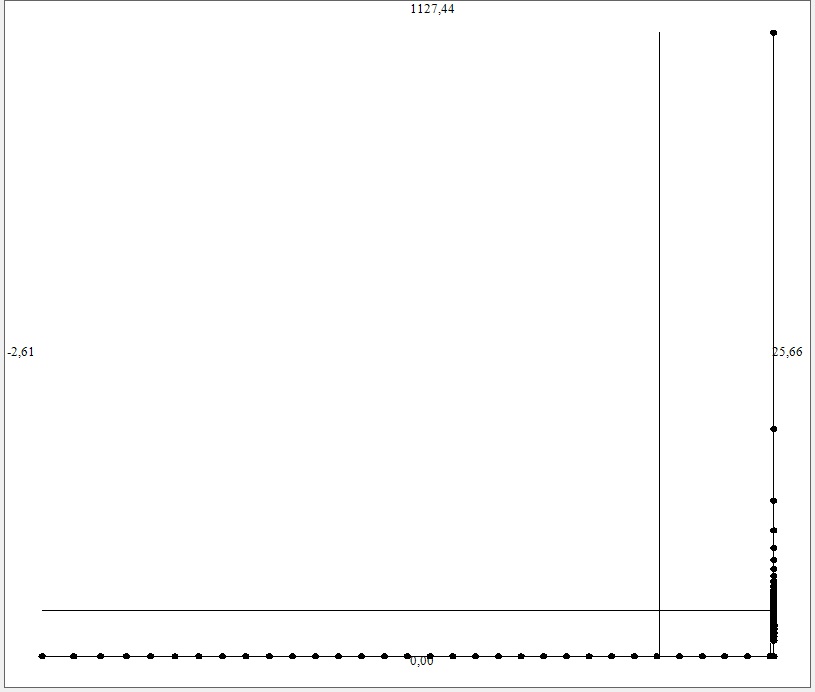}
{(d) t=106}
\end{minipage}

\end{tabular}
\caption{Solution component $\rho(t,\cdot)$} \label{Eulerrho}
\end{figure}

Figures \ref{Euleru} and \ref{Eulerrho} display the results of our numerical method for this problem. We have plot both
 $w(t,\cdot)$ and $\rho(t,\cdot)$, where the initial density is taken as 
\be
\rho(0,x) = {1 \over { \sqrt{ 2 \pi }}} \exp\Big( - {x^2 \over 2 } \Big).
\ee
 The velocity $u(t,\cdot)$ follows Burgers equation, and is computed as we explained earlier.
 (See Figure \ref{BEE} for the velocity at the corresponding times.)
 Note that the density $\rho(t,\cdot)$ is singular after the shock formation, and we thus performed a suitable 
truncation in order to display the function $\rho(t,\cdot)$ in Figure \ref{Eulerrho}.



\subsection{Conservation laws with one inflection point} 

We next investigate the class of conservation laws of the form \eqref{CL} with non-convex flux. 
Such laws arise in, for instance, material science or the dynamics complex fluid flows.

\begin{figure}[h]
\centering
\begin{tabular}{ccc}

\begin{minipage}{1.65in}
\includegraphics[height=1.6in,width=1.9in]
{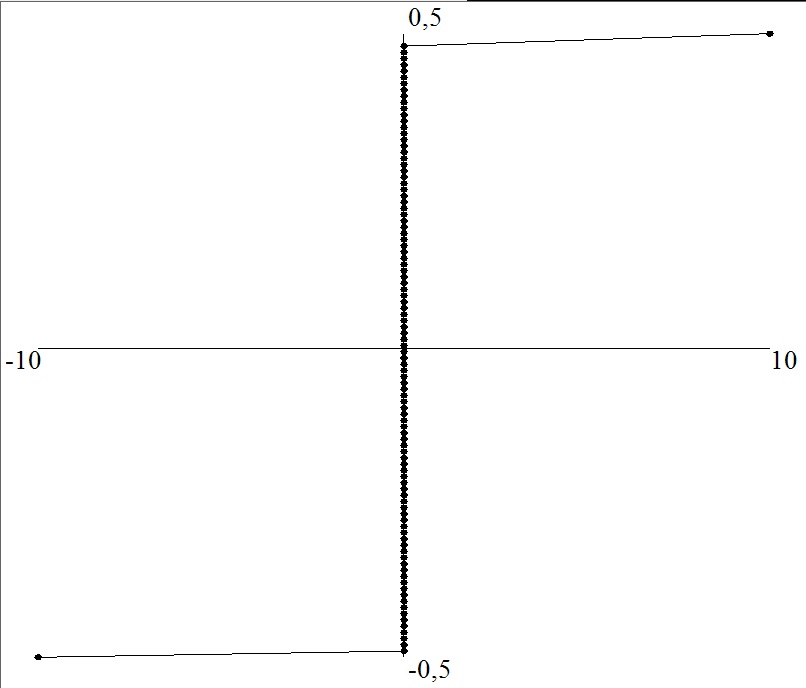}
t=0.
\end{minipage}

\hskip3.01cm

\begin{minipage}{1.6in}
\includegraphics[height=1.6in,width=1.9in]
{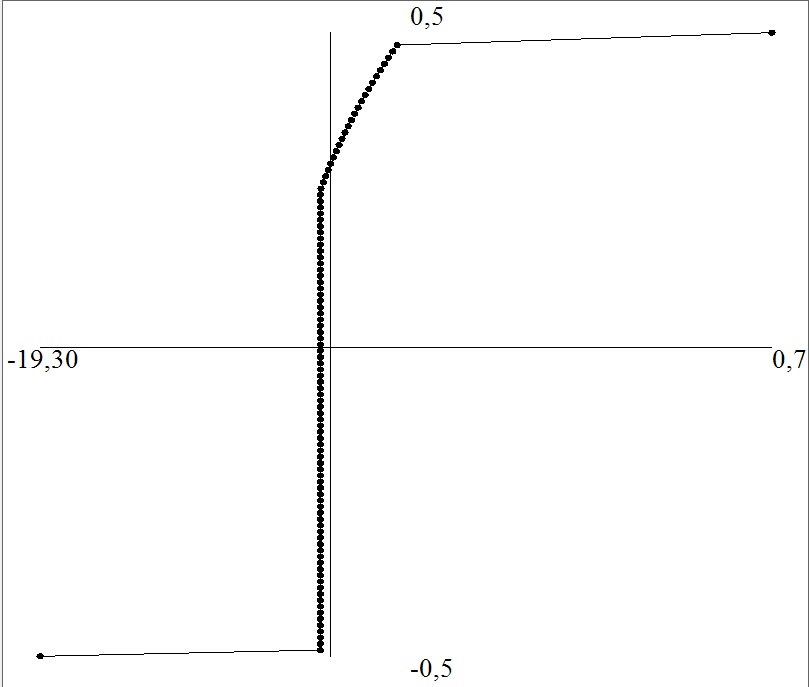}
t=10.
\end{minipage}

\end{tabular}
\caption{Solutions with non-convex flux} \label{NMSI}
\end{figure}

Specifically, in \eqref{eq501}, 
we choose  $f(u) = {1 \over 3}  u^3 - u$. The first numerical test corresponds to Riemann data generating a non-monotone shock formation, that is a shock  $(u_l,u_m = .25)$ followed by a rarefaction wave $(u_m = .25,u_r)$; see Figure \ref{NMSI}: 
$$
	 u(0,x) = u_l := -0.5 \text{ for } x < 0;\qquad   u_r := 0.5 \text{ for } x \ge 0.
$$

\begin{figure}[h]
\centering
\begin{tabular}{ccc}

\begin{minipage}{1.65in}
\includegraphics[height=1.6in,width=1.9in]
{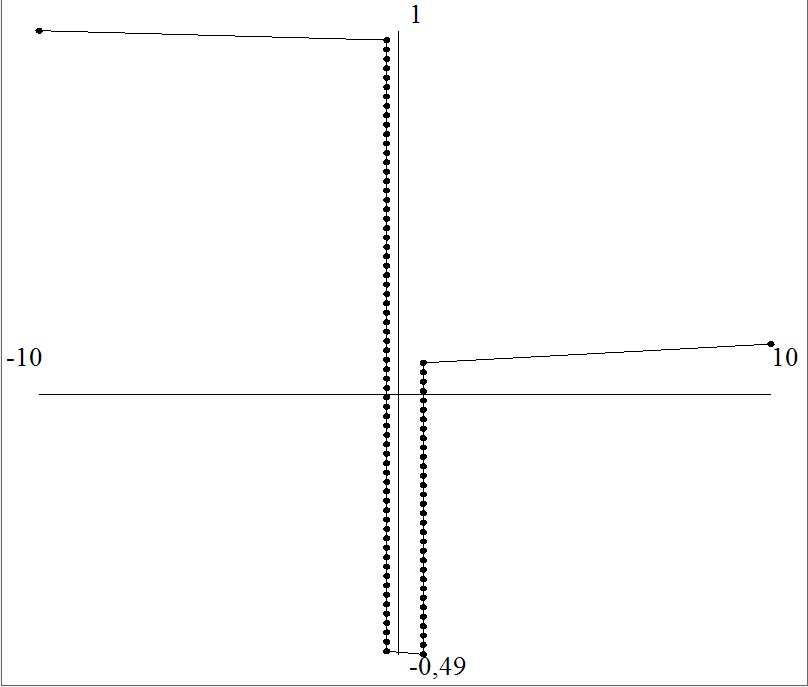}
(a) t=0
\end{minipage}

\hskip3.01cm

\begin{minipage}{1.6in}
\includegraphics[height=1.6in,width=1.9in]
{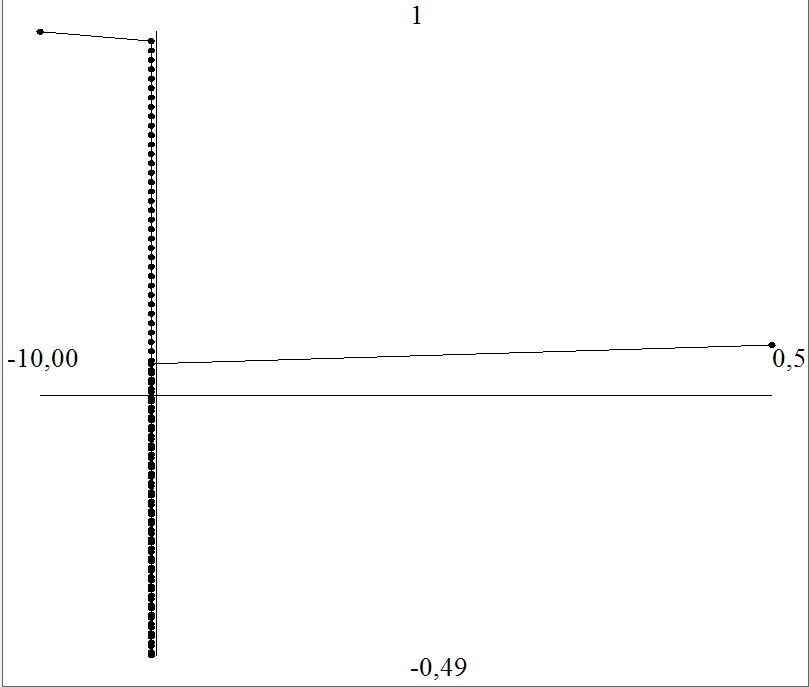}
(b) t=10
\end{minipage}

\end{tabular}
\caption{Shock wave interaction} \label{MSI}
\end{figure}

For our second numerical test, we use an initial data with {\it two} shocks, which initially separated but interact at a later time. After the  interaction, the middle state $u_m$ cancels out and the solution stabilizes to  single shock $(u_l = 1.,u_r := .25)$ traveling at a negative velocity (cf.~Figure \ref{MSI}): 
$$
	 u(0,x) = u_l := 1. \text{ for } x <- 0.5;\quad  u_m := -0.5 \text{ for }   -0.5 \le x \le 0.5;
$$
$$
	 u(0,x) = u_r := .25 \text{ for } x \ge 0.5.
$$
It is remarkable that, even for such a non-convex flux,  our algorithm generates the entropy solution (i.e. physically meaningful) to the problem, which is the one characterized  in \cite{Kruzkov}.  

\subsection{Conservation laws with several inflection points}

\begin{figure}[h]
\centering
\begin{tabular}{ccc}

\begin{minipage}{1.65in}
\includegraphics[height=1.6in,width=1.9in]
{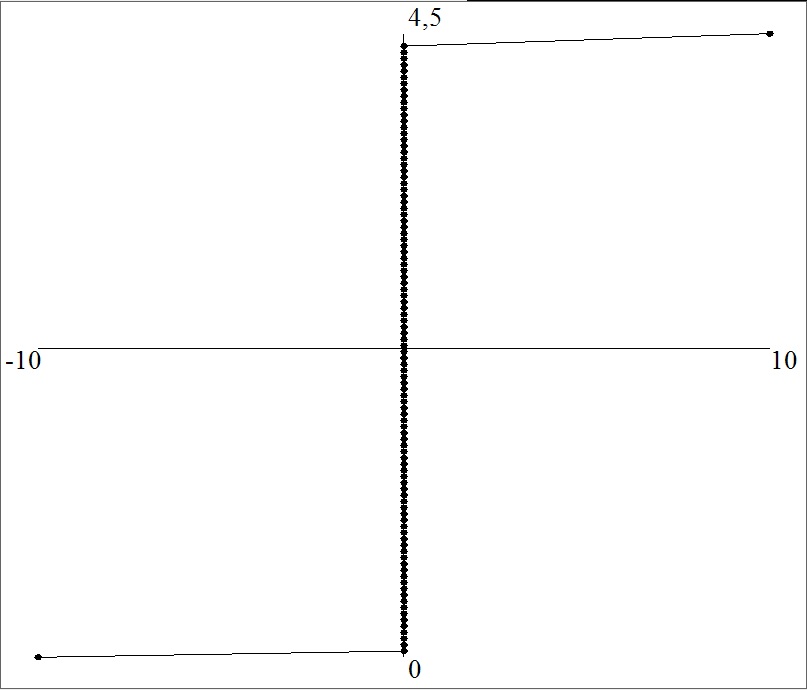}
t=0
\end{minipage}

\hskip3.01cm

\begin{minipage}{1.6in}
\includegraphics[height=1.6in,width=1.9in]
{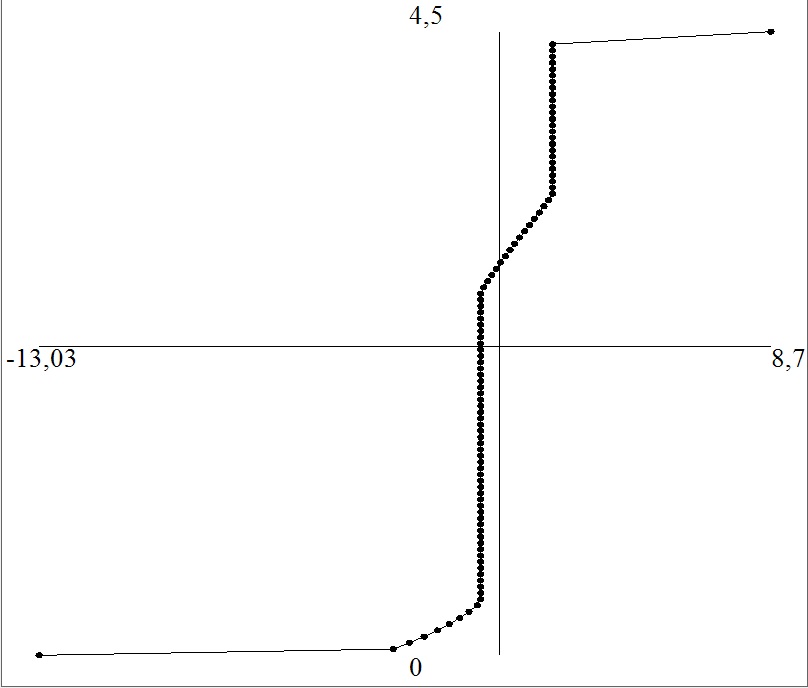}
t=20
\end{minipage}

\end{tabular}
\caption{Complex wave pattern for a flux having multiple inflection points} \label{NMSI_complex}
\end{figure}

We finally apply the CHA method to a conservation law whose flux admits several inflection points. Specifically, we consider the nonconvex function 
\be 
f(u) = -u(u-1)(u-2)(u-3)(u-5) / 200.
\ee
This example provides us with a further challenge for our method. The first test we consider corresponds to the Riemann data 
\be
	 u(0,x) = u_l = 0 \text{ for }  x < 0; \qquad   u_r := 4.5 \text{ for } x > 0, 
\ee
with which the solution is expected to exhibit a rather complex wave pattern, that is, 
 two rarefaction waves and two shocks. The numerical results are strikingly in accordance with the theoretically expected ones; see Figure \ref{NMSI_complex}.

\begin{figure}[h]
\centering
\begin{tabular}{ccc}

\begin{minipage}{1.3in}
\includegraphics[height=1.6in,width=1.9in]
{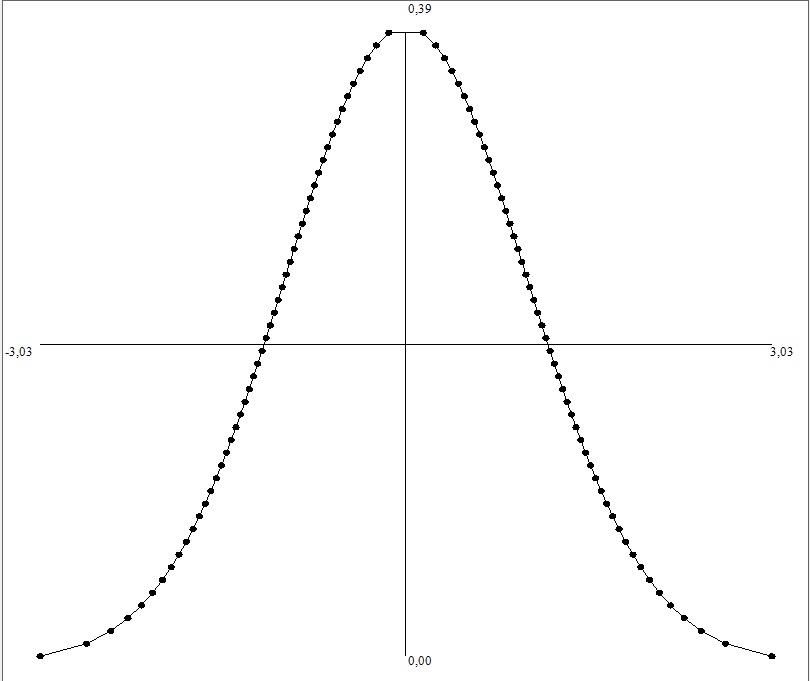}
{t=0.}
\end{minipage}

\hskip3.01cm

\begin{minipage}{1.3in}
\includegraphics[height=1.6in,width=1.9in]
{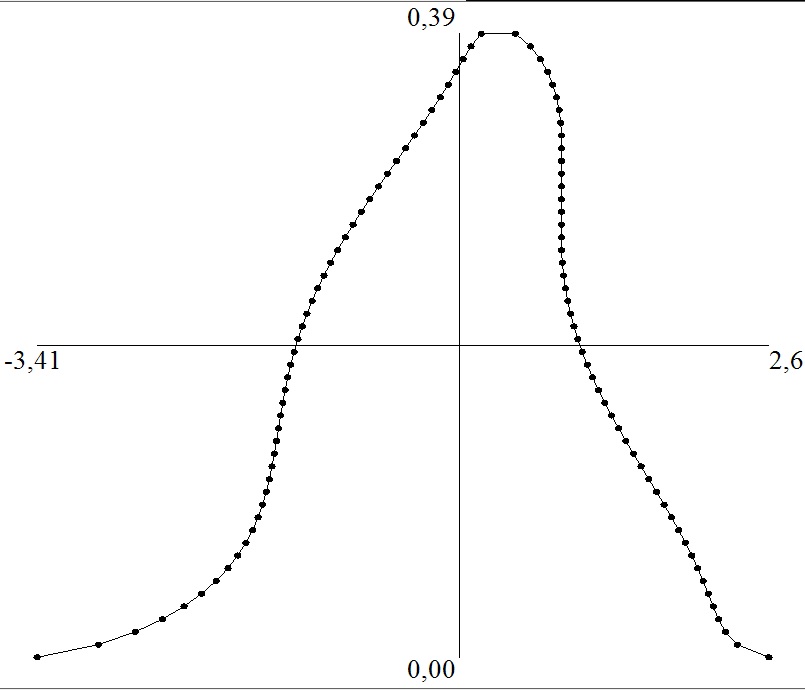}
{t=.5}
\end{minipage}

\\

\begin{minipage}{1.3in}
\includegraphics[height=1.6in,width=1.9in]
{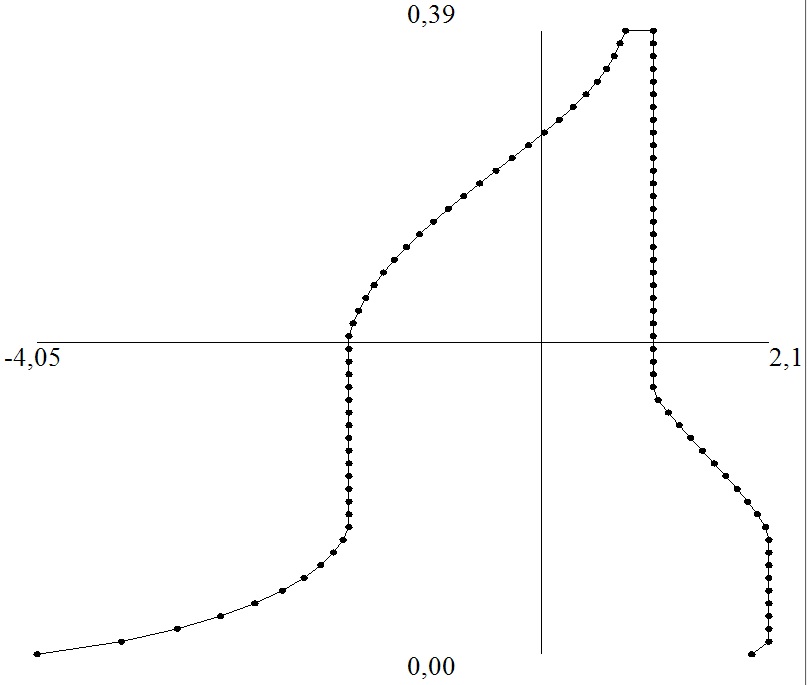}
{t=1.}
\end{minipage}

\hskip3.01cm

\begin{minipage}{1.3in}
\includegraphics[height=1.6in,width=1.9in]
{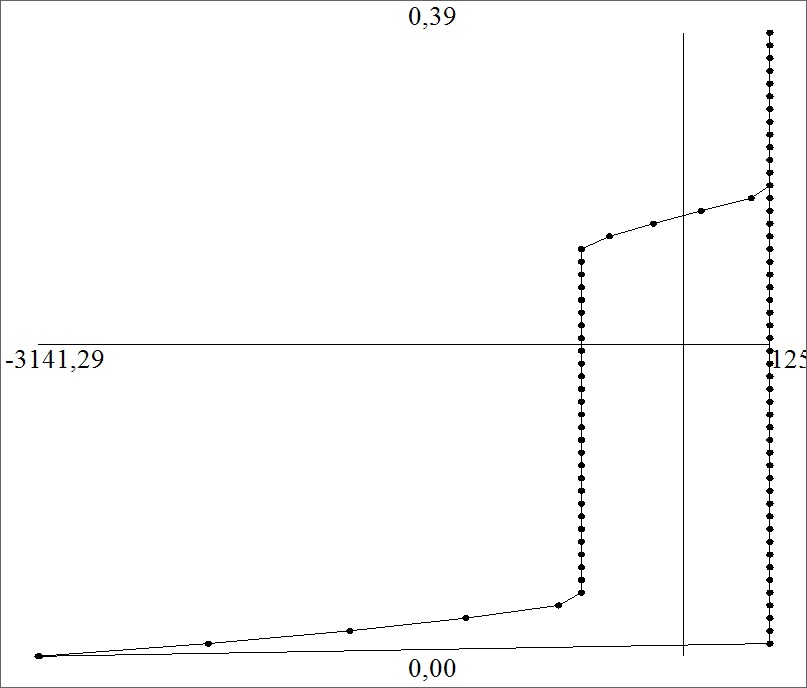}
{t=1001}
\end{minipage}

\end{tabular}
\caption{Asymptotic convergence with a flux having multiple inflection points} \label{u_complex_gaussian_CHA}
\end{figure}

Second, we investigate the time-asymptotic behavior of the CHA-algorithm with the same flux of degre five, but now with a Gaussian initial data (as defined earlier). Even though this data generates a complex wave pattern 
for sufficiently small times, the solutions are expected to converge toward a generalization of the N-wave of
 the convex case. Such a kind of convergence is indeed observed with our method, as
 illustrated by Figure \eqref{u_complex_gaussian_CHA}. Yet, we emphasize the non-trivial shape of the asymptotic N-wave in this case. 

\section{Concluding remarks}

The method proposed in this paper has the following advantages: 

\begin{itemize}

\item The method applies to hyperbolic problems and Hamilton-Jacobi equations. 
It generates solutions at any arbitrary time and does not require a time-evolution. In particular, we can easily compute the solutions at an arbitrary large time and therefore 
determine the asymptotic regime of the solutions (such as the $N$-wave profile for convex hyperbolic equations).

\item The key observation made in this work is that the multivalued solutions contains all the "information" about the entropy dissipative solutions, so that the latter can be recovered from the former by our convex hull algorithm  (CHA).

\item The jump discontinuities in solutions (shock waves) are sharply represented 
and described with finitely many points proportionally to the discontinuity strength. This is in contrast with methods based on finite difference schemes for which the numerical solutions contain only a few points within the shocks.

\item Especially in a multidimensional setting, the sampling  ${\it Y} := \big(Y_n \big)_{n=1, \ldots, N}$ (cf.~Section \ref{algo}) should be chosen to be 
equidistributed in $\Lambda$. Such a set can be computed by a random generator, but a better choice is an optimal quantizer of the "uniform law", which can be achieved when the initial condition is written as a discrete sum of convex and concave components.

\item Optimization strategies can be developed for the proposed algorithm. For instance, removing discrete points lying inside the convex hull and adding discretization points outside. Such a strategy allows for instance to 
keep the discretization error bounded uniformly  in time.
   
\end{itemize} 


\section*{Acknowledgments}

The first author (PLF) was partially supported by the Centre National de la Recherche Scientifique (CNRS) 
and the Agence Nationale de la Recherche through the grants ANR~2006--2-134423 and SIMI-1-003-01.


\end{document}